\documentclass[12pt,twoside]{article}
\usepackage[dvips]{graphicx}
\usepackage{latexsym}

\newtheorem{theorem}{Theorem}
\newtheorem{cor}[theorem]{Corollary}

\newtheorem{lemma}[theorem]{Lemma}

\newtheorem{definition}[theorem]{Definition}

\newtheorem{example}[theorem]{Example}

\begin{document}
\overfullrule=0pt
\baselineskip=24pt
\font\tfont= cmbx10 scaled \magstep3
\font\sfont= cmbx10 scaled \magstep2
\font\afont= cmcsc10 scaled \magstep2
\title{\tfont A Two-Dimensional Minkowski ?(x) Function}
\bigskip
\author{Olga R. Beaver\\ Thomas Garrity\
\\ Department of Mathematics\\ Williams College\\ Williamstown, MA  01267\\
obeaver@williams.edu \\
tgarrity@williams.edu}

\date{}
\maketitle
\begin{abstract}
A one-to-one continuous function from a triangle to itself is defined 
that has both interesting number theoretic and analytic properties.  
This  function is shown to be a natural generalization of the 
classical Minkowski $?(x)$ function.  It is shown there exists a natural class of 
pairs of cubic irrational numbers in the same cubic number field that are 
mapped to pairs of rational numbers, in analog to $?(x)$ mapping 
quadratic irrationals on the unit  interval to rational numbers on 
the unit interval.  It is also shown that this new function satisfies 
an analog to the fact that $?(x)$, while increasing and continuous, 
has derivative zero almost everywhere.

\end{abstract}

\section{Introduction}
Any real number with an eventually periodic continued fraction expansion 
must be a quadratic irrational.
 This property linking periodicity of a number's continued fraction 
 expansion with its being quadratic  led Minkowski 
 to  define his remarkable question-mark function
$$?:[0,1]\rightarrow [0,1].$$
\noindent (See page 50 of
Volume II in \cite{Minkowski1}; see also p. 754, article 196 in 
\cite{Hancock1}, which appears to be essentially a translation of all of 
Minkowski's number theory papers.)  The question-mark function is increasing, continuous, maps 
 each rational number $\frac{p}{q}$
 to a pure dyadic number of the form $\frac{k}{2^{n}},$ maps
each quadratic irrational to a rational number, and  has the property 
that
 the inverse image of the rational numbers is exactly the set of quadratic
 irrationals.  In order to understand the number theoretic properties 
 of quadratic irrationals, it is natural to look at the
 function theoretic properties of $?(x)$.  In particular, the 
 question-mark function is not only continuous and monotonically 
 increasing but 
has derivative zero almost everywhere.  As such, it is a naturally
occuring example of a singular function.   Moreover, it is, in 
fact, the diophantine properties of 
continued fractions
that lead to its derivative being zero a.e.. Thus the analytic 
property of $?(x)$ being both increasing and having derivative zero 
almost everywhere is actually number theoretic in origin.  
 
 In this paper, we will construct 
 a function similar to Minkowski's question-mark function, and will 
 use that function in order to understand the properties of cubic irrationals.

A. Denjoy, in \cite{Denjoy1}, \cite{Denjoy2} and independently R. Salem, in \cite{Salem1},
were the first to realize that $?(x)$ is
singular, although earlier, F. Ryde
\cite{Ryde1} proved in essence that  $?(x)$ was singular.  However, Ryde showed 
that $?(x)$ was singular without realizing its connection
with Minkowski's function (see also Ryde's \cite{Ryde2}). Recent work
on $?(x)$ is the work of Kinney \cite{Kinney1}, Girgensohn \cite{Girgensohn1},
Ramharter \cite{Ramharter1}, of Tichy and Uitz \cite{Tichy-Uitz1},  and of Viader,
Paradis and Bibiloni \cite{Viader-Paradis-Bibiloni1}
\cite{Paradis-Viader-Bibiloni1}. (In fact, the
idea for the inequality that we prove in section 6.1 and use in 6.2
was inspired by the work of Viader,
Paradis and Bibiloni in \cite{Viader-Paradis-Bibiloni1} .)

A natural question to ask is: do cubic irrationals and other higher order 
algebraic numbers lend themselves to similar analysis?  An even more 
basic question to ask is how to generalize the relation between 
periodicity for continued fractions and  quadratic irrationals to 
cubics.  In
1848, in a letter to Jacobi, Hermite \cite{Hermite1} asked for such a 
generalization.
  Specifically,  the Hermite problem 
is:

\noindent {\it Find
methods for expressing real numbers as sequences of positive integers
so that the sequence is eventually periodic precisely when the initial
number is a cubic irrational.}

Over the years there has been much work in trying to solve the Hermite 
problem.  For an overview, see Schweiger's {\it Multidimensional 
Continued Fractions} \cite{Schweiger4}. For work up to 1980, see 
Brentjes' overview in \cite{Brentjes1}.  Other work is in 
\cite{ACCCDGLS1},
\cite{Bernstein1}, \cite{Ferguson-Forcade1},\cite{FHZ1}, \cite{Garrity1}, 
\cite{Grabiner-Lagarias1}, \cite{GM1}, 
\cite{Klein1},\cite{Korkina1},\cite{Korkina2},\cite{Korkina3}, 
\cite{Lachaud1}, \cite{Lagarias1}, \cite{Mahler1}, 
\cite{Schweiger1},\cite{Schweiger2},\cite{Schweiger3}.

  On the other hand, there has been little attempt to 
approach the Hermite problem by 
generalizing the Minkowski $?(x)$.  The only such attempt that we have
found is in the thesis of Louis Kollros \cite{Kollros1}.  Kollros
generalizes $?(x)$ to a map from the unit square to itself.  However, 
while he sets up various methods for associating points in
the unit square with sequences of integers, 
he does not concern himself with the function-theoretic properties of this
function.  It does not appear that
Kollros has solved the Hermite problem.  In particular, he was not interested
in the differentiability properties of his analogue to $?(x)$.

In this paper we develop a different, more natural, analog to $?(x)$.  
In section two, a review of the Minkowski question-mark function is 
given.  In section three,  we
 construct a
map from a two dimensional simplex (a triangle) to itself, as an analog
to the map of $?(x)$ from a one dimensional simplex to itself.  The
map will be determined by partitioning the triangle, first via a ``Farey''
partition, and then by a barycentric (triadic) partitioning, which we 
will frequently call the ``bary'' partitioning.  We define 
a function $\delta(x,y)$ from the Farey triangle to the barycentric 
triangle.  In section four, we see that the Farey iteration can be 
viewed as a multidimensional continued fraction.  We show that 
periodicity of the Farey iteration corresponds to a class of cubic 
irrationals. In section five we show, by contrast, that periodicity 
for the barycentric iterations corresponds to a class of rational 
points.  This results in  that our function will  map a natural class of cubic 
points 
to a natural class of rational points.  Finally, in section six, we  prove an analog of
singularity by showing that, a.e., the area of image triangles in the
barycentric partitioning approaches zero far more quickly than the
area of the
domain triangles in the Farey partitioning.

We note that using Farey partitioning, or Farey nets, to solve the Hermite problem has been 
considered by both Monkemeyer \cite{Monkemeyer1} and
more recently by Grabiner \cite{Grabiner1}. Both papers are quite
interesting; neither use Farey nets to generalize
the Minkowski $?(x)$ function.  In actuality, this analytic approach would not have been a 
natural succession in either of these papers, as Monkemeyer's and Grabiner's goals were not function
theoretic.

We would like to thank Lori Pedersen for making all of the diagrams.  
Also, we would like to thank Keith Briggs for pointing out some 
errors in the bibliography of an earlier version of this paper.

\section{A Review of the Minkowski Question-Mark Function}

All of the discussion in this section is well-known.  We include it 
here for sake of completeness.  

Recall that given two rational numbers $\frac{p_{1}}{q_{1}}$ and
$\frac{p_{2}}{q_{2}}$, each in lowest terms, the {\it Farey sum}, $\hat{+}$, of the numbers is
$$\frac{p_{1}}{q_{1}} \hat{+}\frac{p_{2}}{q_{2}} =
\frac{p_{1}+p_{2}}{q_{1}+q_{2}}.$$
The ? function is then defined as follows.  Suppose we know the value of $?(\frac{p_{1}}{q_{1}})$ and
$?(\frac{p_{2}}{q_{2}})$.  We then set
$$?(\frac{p_{1}+p_{2}}{q_{1}+q_{2}})= \frac{?(\frac{p_{1}}{q_{1}})+
?(\frac{p_{2}}{q_{2}})}{2}.$$
Specifying the intial values
$$?(0)=0\;\mbox{and}\; ?(1)=1,$$
we now know the values of $?(x)$ for any rational number $x$.

 By
continuity arguments we can determine the values of $?(x)$ for any
real number $x$ in the unit interval.  
Since we will be generalizing
this continuity argument in the next section, we discuss this now in 
some
 detail.

We produce two sequences of partitions, $\mathcal{I}_{k}$ and
$\tilde{\mathcal{I}}_{k}$, of the unit interval. For each $k\geq 0$, each
partition will split the unit interval into $2^{k}$ subintervals. Both
start with
just the unit interval itself:
$$\mathcal{I}_{0}=\tilde{\mathcal{I}}_{0}=[0,1].$$  Note that
$0=\frac{0}{1}$ and $1=\frac{1}{1}$.  Now, given the partition $\mathcal{I}_{k-1}$, suppose that the endpoints
of each of its $2^{k-1}$ open subintervals are rational numbers.  Form the next partition $\mathcal{I}_{k}$ by taking the Farey sum of
the endpoints of the partition $\mathcal{I}_{k-1}$.  Thus the
endpoints of $\mathcal{I}_{k}$  consist of the Farey fractions of
order $k$.

\begin{tabbing}
\bigskip
\hspace*{0.375in}\=
\scalebox{0.55}
{\includegraphics*[bb=0in 3.75in 8.5in 9.in]{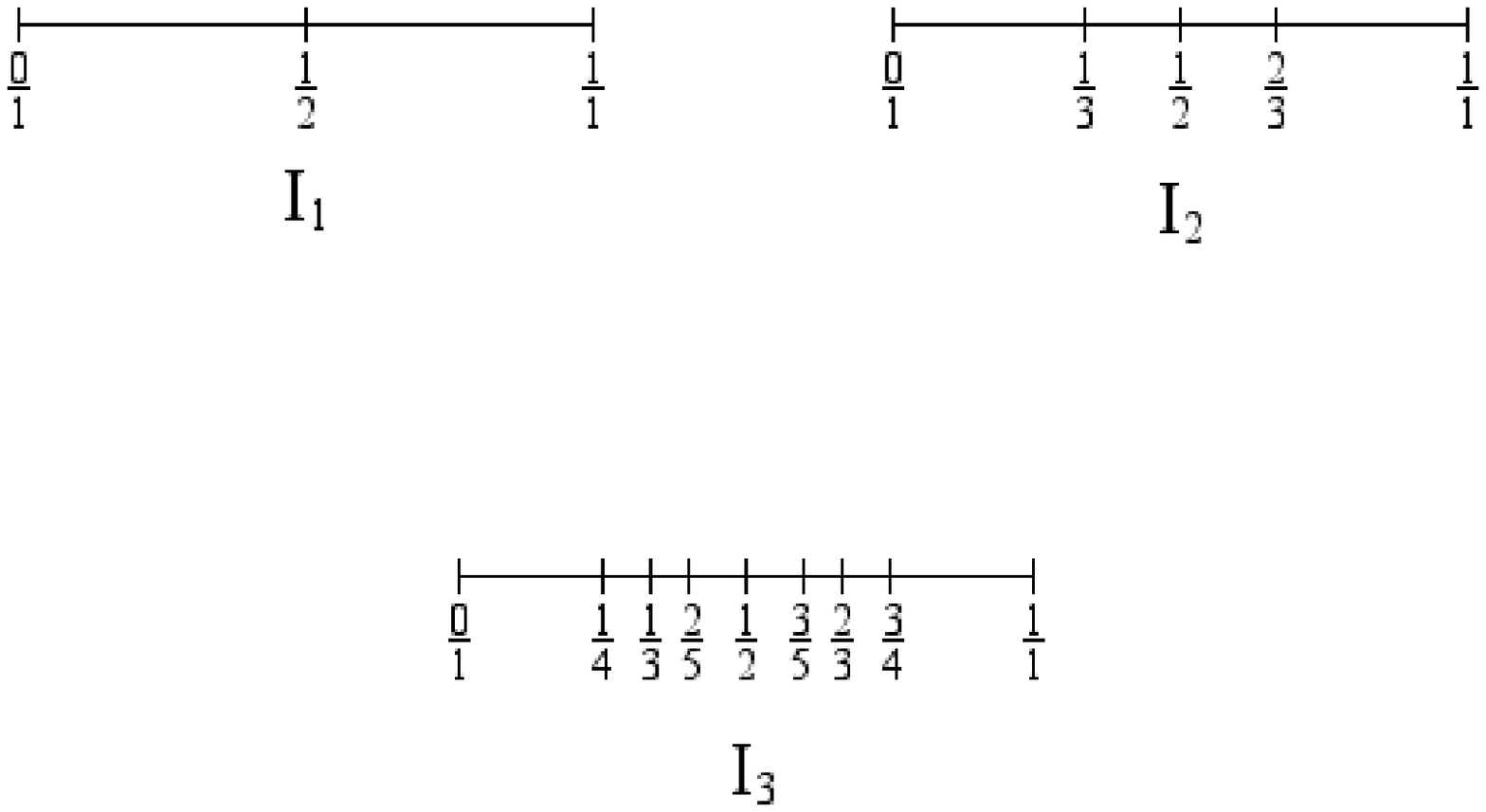}}
\bigskip
\end{tabbing}

The partition $\tilde{\mathcal{I}}_{k}$ is even simpler.  It is just the partition
given by the subintervals $[\frac{l-1}{2^{k}},\frac{l}{2^{k}}]$.

\begin{tabbing}
\bigskip
\hspace*{0.375in}\=
\scalebox{0.55}
{\includegraphics*[bb=0in 4.25in 8.5in 9.25in]{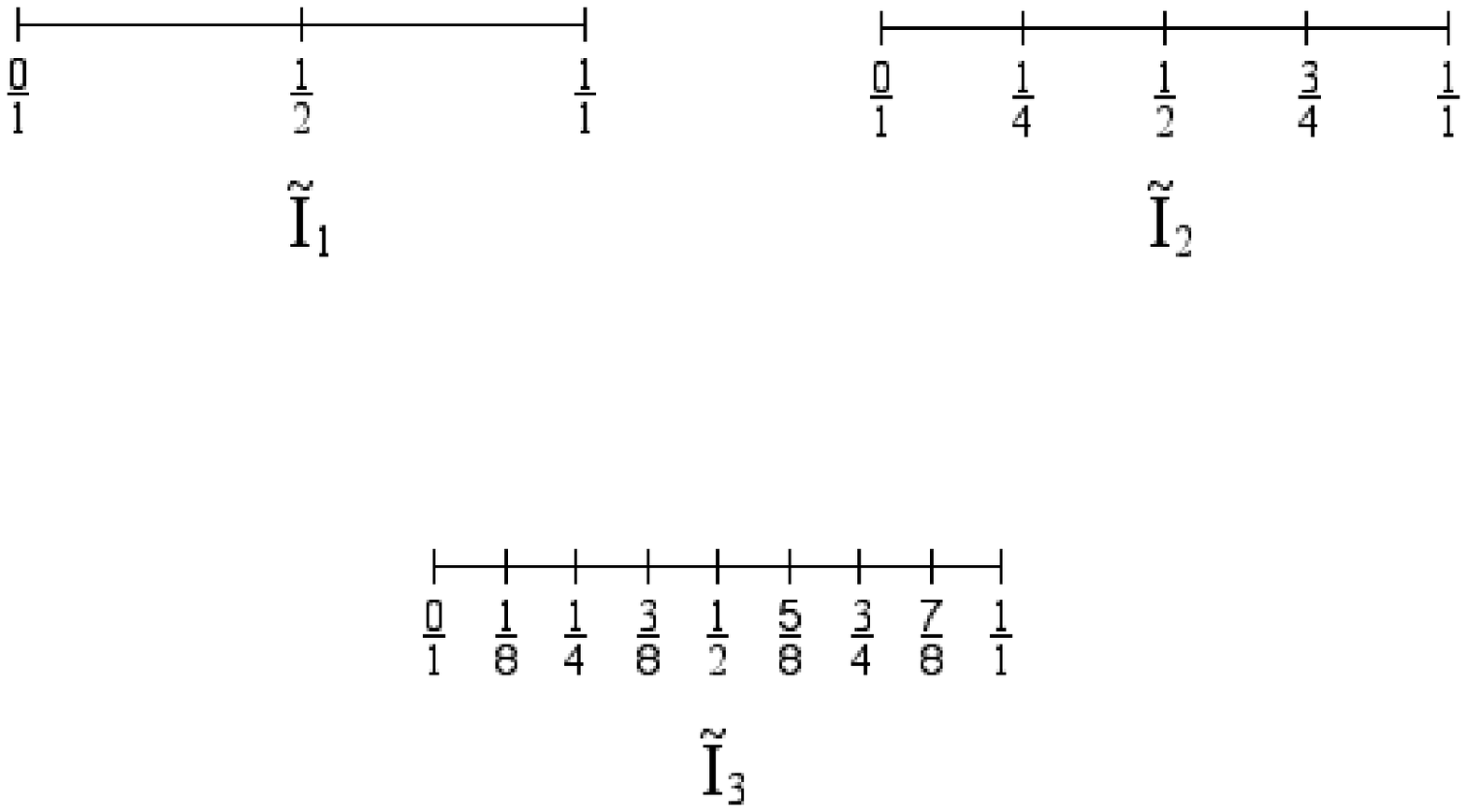}}
\bigskip
\end{tabbing}

Then the function $?(x)$ can be seen to map the endpoints of each 
$\mathcal{I}_{k}$ to the corresponding endpoints of $\tilde{\mathcal{I}}_{k}$.

Now, as is shown, for example, in \cite{Salem1},  $?(x)$ is singular and
hence has derivative zero almost everywhere.  Using the partitions 
defined above, we can recast the fact that   $?(x)$ is a singular 
function into the language of lengths of intervals.
Fix $\alpha\in [0,1]$.  For each $k$, let $I_{k}$ and $\tilde{I}_{k}$
be the subintervals of the respective partitions $\mathcal{I}_{k}$ and
$\tilde{\mathcal{I}}_{k}$ that contain the point $\alpha$.
Then, as shown page 437 in \cite{Salem1},
\begin{theorem}
For almost all $\alpha\in [0,1]$,
$$\liminf_{k\rightarrow \infty} \frac{\mbox{length of}\; \tilde{I}_{k}}
{\mbox{length of}\;I_{k}} = 0.$$
\end{theorem}
It is this theorem that provides the most natural language for 
generalizing the failure of differentiability for our analog of the question-mark function.

The proof involves the idea that the Diophantine approximations 
properties of continued fractions make the above denominator approach 
zero more slowly than the numerator.

 
\section{The Farey-Bary Map: A generalization of the Minkowski 
Question-Mark Function.}  

Our goal is to define a continuous map from a two-dimensional 
simplex (a triangle) to itself that generalizes 
the 
Minkowski question-mark function.  This will involve two separate 
partitionings of the triangle.  We would like to have 
  periodicity in 
the domain  correspond to  cubic irrationals while periodicity in the range 
to imply rationality.  Both of these goals will only be achieved in 
part, as we will show that periodicity will imply cubic irrationality 
in the domain case 
and rationality for the range. At the same time, we want our generalization to 
obey some sort of singularity property. 

Although in a sense it would be most natural to denote our generalization
by the symbol $?(x,y)$, we have found that it is both awkward to say and
awkward to read.  Thus we will denote our generalization by $\delta (x,y)$.


\subsection{The Farey Sum in the Plane}

We will often need to refer to a point in the plane of the form 
$${v} = \pmatrix{p/r\cr 
q/r\cr }.$$  Here, since the coordinates share the same denominator, we 
can associate to this point   a unique  vector in space, namely  
$$\bar{v}= 
\pmatrix{p\cr q\cr r\cr}.$$  Conversely, a vector  $\bar{v}= 
\pmatrix{p\cr q\cr r\cr}$ can be associated uniquely to the point  
$${v} = \pmatrix{p/r\cr q/r\cr }$$ in the plane.   In what follows, we will usually refer to 
both the point and its corresponding vector as $v$.

Consider three points in the plane, each of whose entries are nonnegative 
integers, each $r_{i}\neq 0$,
and such that each vector's entries share no common factors:
$${v_{1}} = \pmatrix{p_{1}/r_{1}\cr q_{1}/r_{1}\cr },
{v_{2}}= \pmatrix{p_{2}/r_{2}\cr q_{2}/r_{2}\cr },
{v_{3}}= \pmatrix{p_{3}/r_{3}\cr q_{3}/r_{3}\cr }.$$  These points 
define a triangle in the plane and, as noted above, can also be represented as the vectors, 
$$v_{1}= \pmatrix{p_{1}\cr q_{1}\cr r_{1}\cr}, v_{2}=\pmatrix{p_{2}\cr
q_{2}\cr r_{2}\cr} \mbox{and }v=\pmatrix{p_{3}\cr q_{3}\cr 
r_{3}\cr}.$$  Summing the three vectors, we get
$$v = v_{1} + v_{2} + v_{3} =
\pmatrix{p_{1}+ p_{2}+p_{3}\cr q_{1}+q_{2}+q_{3}\cr 
r_{1}+r_{2}+r_{3}\cr}.$$  This vector sum can be converted into a 
point v in the plane, where
$$v = \pmatrix{\frac{p_{1}+p_{2}+p_{3}}{r_{1}+r_{2}+r_{3}}\cr \cr 
\frac{q_{1}+q_{2}+q_{3}}{r_{1}+r_{2}+r_{3}} }.$$ 
This correspondence between points in the 
plane and  a vector representation allows us to define the Farey sum.

\begin{definition}
Let  
$$v_{1} = \pmatrix{p_{1}/r_{1}\cr q_{1}/r_{1}\cr },
{v_{2}}= \pmatrix{p_{2}/r_{2}\cr q_{2}/r_{2}\cr },
{v_{3}}= \pmatrix{p_{3}/r_{3}\cr q_{3}/r_{3}\cr },$$
where, for each $i$, the $p_{i}$, $q_{i}$ and $r_{i}$ share no common 
factor. 
The Farey sum, $\hat{v}$ of the  $v_{i}$ is then

$$\hat{v}=v_{1} \hat{+} v_{2} \hat{+} v_{3} =\pmatrix{\frac{p_{1}+ p_{2}+p_{3}}{r_{1}+r_{2}+r_{3}}
\cr \cr \frac{q_{1}+q_{2}+q_{3}}{r_{1}+r_{2}+r_{3}}\cr }.$$

\end{definition}
Note that the point $ \hat{v}$
is inside the triangle determined by the vertices  $v_{1}, v_{2} \mbox{ 
and } v_{3}$ and that $\hat{v}$  corresponds to the vector
 $\pmatrix{p_{1}+ p_{2}+p_{3}\cr q_{1}+q_{2}+q_{3}\cr 
r_{1}+r_{2}+r_{3}\cr}.$


\subsection{Farey and Barycentric Partitions}In this section we will 
define two partitions of the triangle 
$$\bigtriangleup = \{(x,y):1\geq x \geq y \geq 
0\}.$$
 The first partition of $\bigtriangleup$ will yield the domain of our desired 
function $\delta,$ while the second partition will yield the range.

\bigskip

\noindent {\bf The Farey Partition}

\noindent

We will define a sequence of partitions $\{\mathcal{P}_{n}\}$ such 
that each ${\mathcal{P}_{n}}$ will consist of $3^{n}$ subtriangles of 
$\bigtriangleup$ and each ${\mathcal{P}_{n}}$ will be a refinement 
of the previous ${\mathcal{P}_{n-1}}$.

Let ${\mathcal{P}_{0}}$ be the initial triangle $\bigtriangleup$. The three vertices of $\bigtriangleup$ are
$$v_{1}= \pmatrix{0\cr 0\cr }= \pmatrix{0/1\cr 0/1\cr },$$ 
$$v_{2}=\pmatrix{1\cr 0\cr }= \pmatrix{1/1\cr 0/1\cr },$$ $$v_{3}= 
\pmatrix{1\cr 1\cr }=\pmatrix{1/1\cr 1/1\cr }.$$   Taking the Farey 
sum of these vertices,  we have 
 $$v_{1} \hat{+} v_{2} \hat{+} v_{3} =
\pmatrix{2\cr 1\cr 3\cr}.$$  
This Farey vector corresponds to the 
point $\pmatrix{2/3\cr 1/3\cr }.$ In particular, the point 
$\pmatrix{2/3\cr 1/3\cr }$ is an interior point of the triangle 
$\bigtriangleup$ and, in a natural way, partitions $\bigtriangleup$ 
into three subtriangles.  We will refer to the resulting interior 
point as the {\it Farey-center}.

\begin{tabbing}
\bigskip
\hspace*{0.75in}\=
\scalebox{0.5}
{\includegraphics*[bb=0in 6.25in 8.5in 9.9in]{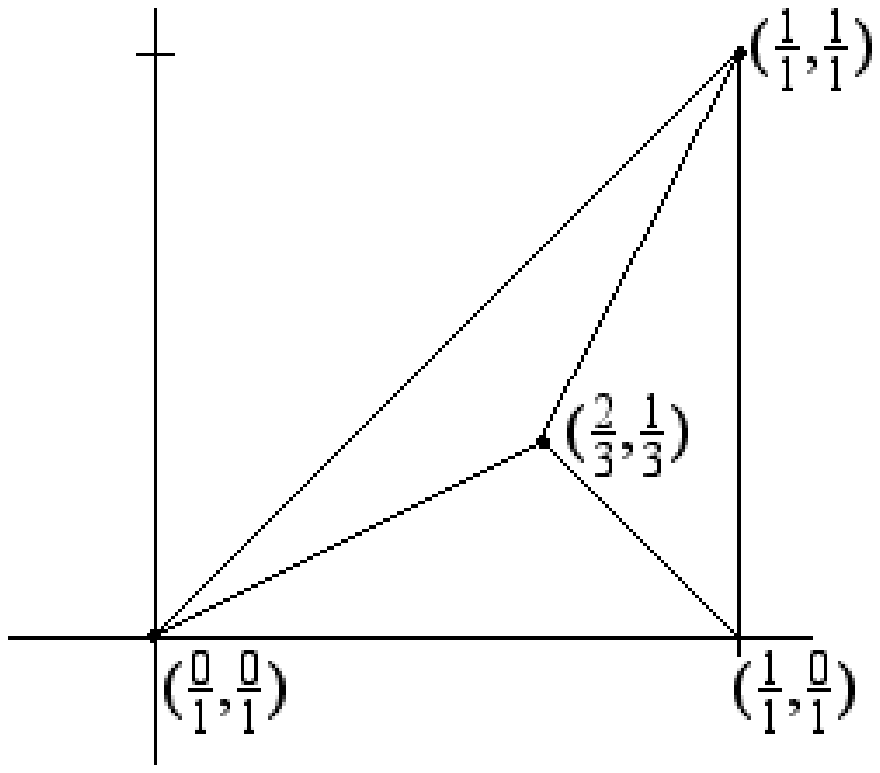}}
\bigskip
\end{tabbing}

This determines the partition $\mathcal{P}_{1}$.  We now proceed 
inductively.  Suppose we have the partition $\mathcal{P}_{n}$, that
determines $3^{n}$ triangles.  We now  partition each of these triangle 
into three subtriangles, as follows.  Suppose one of the triangles in $\mathcal{P}_{n}$
has vertices $\pmatrix{p_{1}/r_{1}\cr q_{1}/r_{1}\cr},
\pmatrix{p_{2}/r_{2}\cr q_{2}/r_{2}\cr }\mbox{and }\pmatrix{p_{3}/r_{3}\cr q_{3}/r_{3}\cr }.$ 
Computing the Farey sum of the three vertices of the 
triangle gives a point $\hat{v}$, the Farey-center, in the interior of 
the subtriangle.  The Farey-center, $\hat{v},$ 
yields a partition of the subtriangle.   Computing in this way the partition of each 
subtriangle of $\bigtriangleup$ determined by $\mathcal{P}_{n},$ gives us the desired next partition 
$\mathcal{P}_{n+1}$ of $\bigtriangleup.$  

We denote this full partitioning of $\bigtriangleup$ by 
$\bigtriangleup_{F}$ and call it the {\it Farey} partitioning.




%


\bigskip

\noindent {\bf The Barycentric Partition}

Again we will define a sequence of partitions 
${\tilde{\mathcal{P}}}_{n}$ of $\bigtriangleup$ such  that each ${\tilde{\mathcal{P}}}_{n}$ will consist of $3^{n}$ triangles of 
$\bigtriangleup$ and each ${\tilde{\mathcal{P}}}_{n}$ will be a refinement 
of the previous ${\tilde{\mathcal{P}}}_{n-1}$.

As with the Farey partitioning, the zeroth level partition 
${\tilde{\mathcal{P}}}_{0}$ is simply the initial triangle 
$\bigtriangleup$.

To compute ${\tilde{\mathcal{P}}}_{1}$, we again start with the original 
three vertices of $\bigtriangleup$ and compute the barycenter of $\bigtriangleup, \mbox{ namely } \pmatrix{2/3\cr 
1/3\cr}$.  This point, called the {\it barycenter}, happens, in this 
case, to be
 the same point obtained by computing the 
Farey sum of the coordinates of the vertices of $\bigtriangleup,$ the 
Farey-center.
This is just a coincidence.

Proceed inductively as follows.  Assume we have a partition 
${\tilde{\mathcal{P}}}_{n}$  of $\bigtriangleup$ into $3^{n}$ 
subtriangles. Further, assume at the $n^{th}$ stage that 
the coordinates of the vertices of any subtriangle can be expressed as rational numbers with 
$3^{n}$ in the denominator. 
 Then, if a given subtriangle in ${\tilde{\mathcal{P}}}_{n}$ has 
vertices $\pmatrix{a_{1}/3^{n}\cr b_{1}/3^{n}\cr }, 
\pmatrix{a_{2}/3^{n}\cr b_{2}/3^{n}\cr },\mbox{ and }
\pmatrix{a_{3}/3^{n}\cr b_{3}/3^{n}\cr },$  we compute the Farey sum of 
the three vertices.  This again gives the barycenter of the subtriangle, 
namely, 
$\pmatrix{\frac{a_{1}+a_{2}+a_{3}}{3^{n+1}}\cr\cr \frac{b_{1}+b_{2}+b_{3}}{3^{n+1}}\cr }.$   Computing, 
in this way, 
the partition of each subtriangle of $\tilde{\bigtriangleup}$ 
determined by ${\tilde{\mathcal{P}}}_{n}$
gives us the desired next partition 
${\tilde{\mathcal{P}}}_{n+1}$ of $\bigtriangleup.$

We call this full partitioning of $\bigtriangleup$ the {\it barycentric}, 
or {\it Bary}, 
partitioning and denote it by $\bigtriangleup_{B}$.

\subsection{The Farey-Bary Map}We are now ready to define the
extension of the Minkowski question-mark function to
$\delta:\bigtriangleup_{F} \rightarrow \bigtriangleup_{B}.$  We will
proceed in stages.  At first, we will  define
a function $\delta_{n}(x,y)$ from the
vertices of the $n^{th}$ partition of $\bigtriangleup_{F}$  to the 
vertices of the $n^{th}$ partition of $\bigtriangleup_{B}$ and then 
extend linearly $\delta_{n}$ to all of $\bigtriangleup_{F}$.
Then
we will show that the functions in the sequence
${\{\delta_{n}}\}$ are continuous and uniformly convergent.  The limit
will be our desired function $\delta(x,y)$ on
$\bigtriangleup_{F}.$

We first need to introduce some notation.  Each of the partitions
$\mathcal{P}_{n}$ and $\tilde{\mathcal{P}_{n}}$ determine subtriangles
of $\bigtriangleup_{F}\mbox{ and } \bigtriangleup_{B},$ respectively.
Let $\bigtriangleup_{n,F}\mbox{ and } \bigtriangleup_{n,B}$ denote
$\bigtriangleup_{F}\mbox
{ and } \bigtriangleup_{B}$ after the n$^{th}$ partitioning,
respectively.  The expression $\langle v_{1}(n), v_{2}(n),
v_{3}(n)\rangle$ will denote a general subtriangle of
$\bigtriangleup_{n,F}$ with vertices $v_{1}(n), v_{2}(n),$ and
$v_{3}(n).$  When we need to refer to the 3$^{n}$ specific
subtriangles, we will use
$\langle{v_{1}}^{s}(n), {v_{2}}^{s}(n),
{v_{3}}^{s}(n)\rangle,$ where s $\in\{1,\ldots,3^{n}\}.$  In a similar
fashion, we will refer to the subtriangles of
$\tilde{\mathcal{P}_{n}}$ by $\langle\tilde{v}_{1}(n), \tilde{v}_{2}(n),
\tilde{v}_{3}(n)\rangle,$ in the general case, and $\langle
{\tilde{v}_{1}}^{s}(n), {\tilde{v}_{2}}^{s}(n),
{\tilde{v}_{3}}^{s}(n)\rangle$ in the specific case.

Note that it happens to be the case that 
$$\mathcal{P}_{0}=\tilde{\mathcal{P}_{0}}, \;
\mathcal{P}_{1}=\tilde{\mathcal{P}_{1}}\;\mbox{and}\;
\bigtriangleup_{1,F}=\bigtriangleup_{1,B}.$$

\begin{definition}Define $\delta_{0},\delta_{1}:\bigtriangleup_{F} \rightarrow 
\bigtriangleup_{B}$ to be the identity maps on the vertices of the 
subtriangles determined by $\mathcal{P}_{0}$ and $\mathcal{P}_{1}$.  
For any $n$, define $\delta_{n}$ to send any vertex in the $n^{th}$ 
partition $\mathcal{P}_{n}$ to the corresponding vertex in the 
partition $\tilde{\mathcal{P}_{n}}$. That is, define 
$\delta_{n}$ on any subtriangle $\langle v_{1}(n), v_{2}(n), 
v_{3}(n)\rangle$ of the partition $\mathcal{P}_{n}$ by  
$$\delta(v_{i}(n))=\tilde{v}_{i}(n)$$ for i = 1,2,3.  Finally, for 
any point $(x,y)$ in the subtriangle with vertices $\langle v_{1}(n), v_{2}(n), 
v_{3}(n)\rangle$, set
$$\delta(x,y)=\alpha\tilde{v}_{1}(n)+\beta\tilde{v}_{2}(n)+\gamma\tilde{v}_{3}(n),$$ 
where
$$(x,y)=\alpha{v}_{1}(n)+\beta{v}_{2}(n)+\gamma{v}_{3}(n).$$
\end{definition}

Note that, since the point $(x,y)$ is in the interior of 
the triangle $\langle v_{1}(n), v_{2}(n), 
v_{3}(n)\rangle$,  we have that 
$$\alpha +\beta + \gamma =1$$
with 
$$0\leq \alpha, \beta, \gamma \leq 1.$$

As defined, $\delta_{0}\mbox{ and }\delta_{1}$ are 
both the identity map since 
$\bigtriangleup_{0,F}=\bigtriangleup_{0,B}$ and 
$\bigtriangleup_{1,F}=\bigtriangleup_{1,B}.$  However, the mappings start to become 
more complicated with $\delta_{2}.$

\pagebreak
\noindent At this stage we have  the 
Farey partition: 

\begin{tabbing}
\bigskip
\hspace*{0.375in}\=
\scalebox{0.5}
{\includegraphics*[bb=0in 3.25in 8.5in 9.75in]{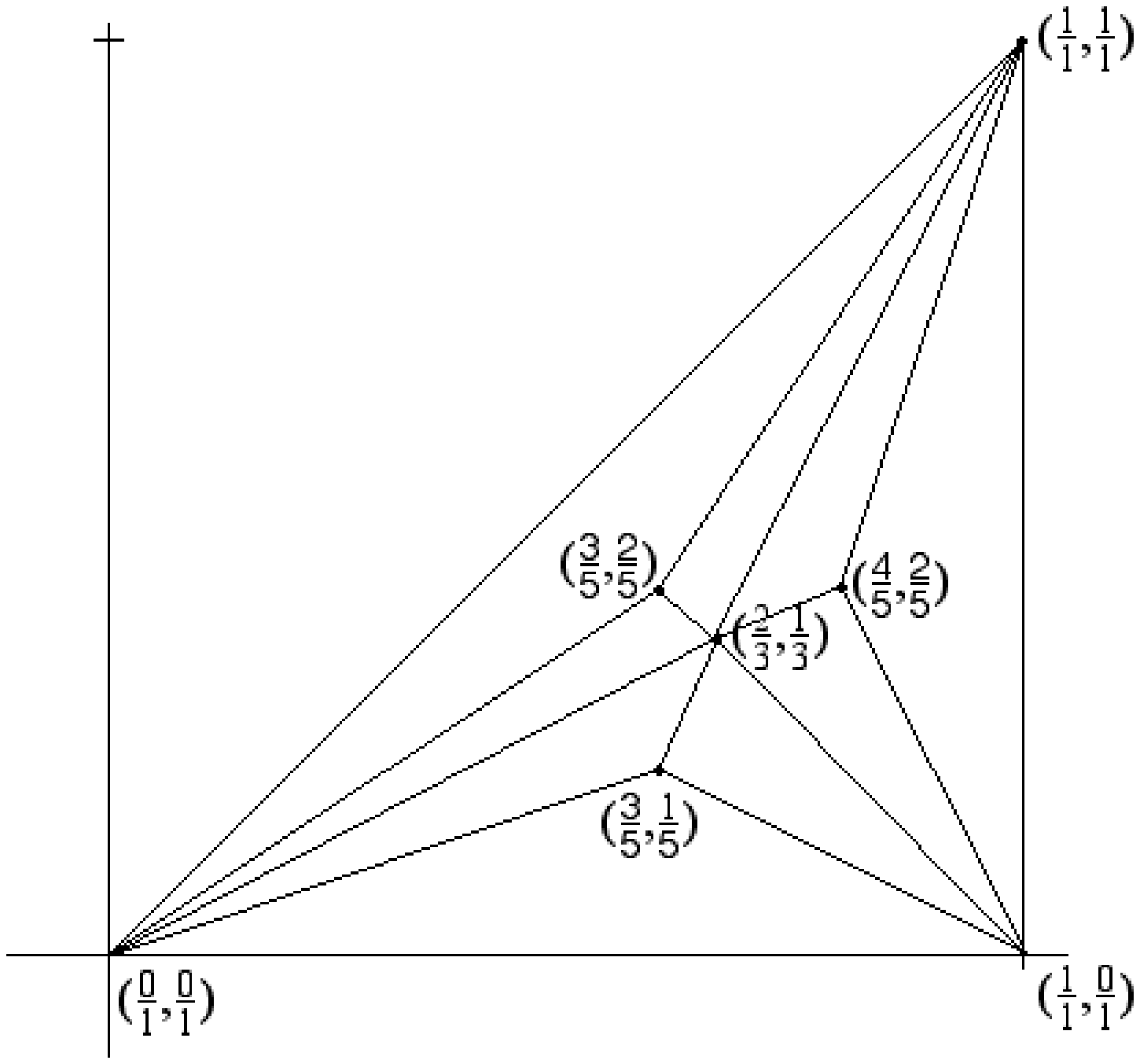}}
\bigskip
\end{tabbing}

\noindent and  the Bary partition:

\begin{tabbing}
\bigskip
\hspace*{0.375in}\=
\scalebox{0.5}
{\includegraphics*[bb=0in 3.25in 8.5in 9.75in]{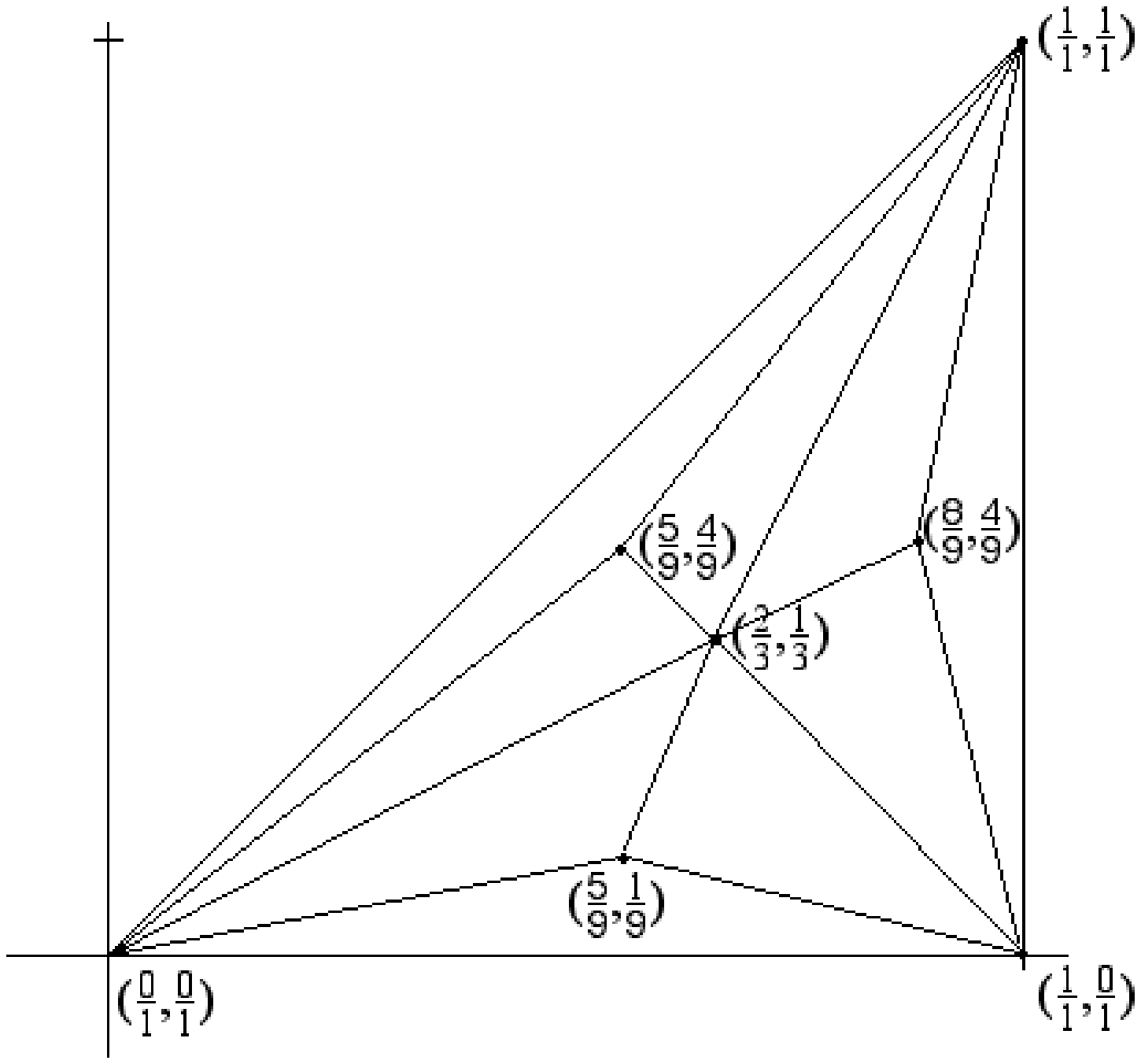}}
\bigskip
\end{tabbing}

\noindent The correspondence between the vertices becomes,\begin{eqnarray*}
\delta_{2} \pmatrix{2/3\cr 1/3\cr } &=& \pmatrix{2/3\cr 1/3\cr }\\
\delta_{2} \pmatrix{3/5\cr 1/5\cr } &=& \pmatrix{5/9\cr 1/9\cr } \\
\delta_{2} \pmatrix{4/5\cr 2/5\cr } &=& \pmatrix{8/9\cr 4/9\cr } \\
\delta_{2} \pmatrix{3/5\cr 2/5\cr } &=& \pmatrix{5/9\cr 4/9\cr }.
\end{eqnarray*} 
Going a few stages further, we get for the Farey partition:

\begin{tabbing}
\bigskip
\hspace*{0.375in}\=
\scalebox{0.5}
{\includegraphics*[bb=0in 3.25in 8.5in 9.75in]{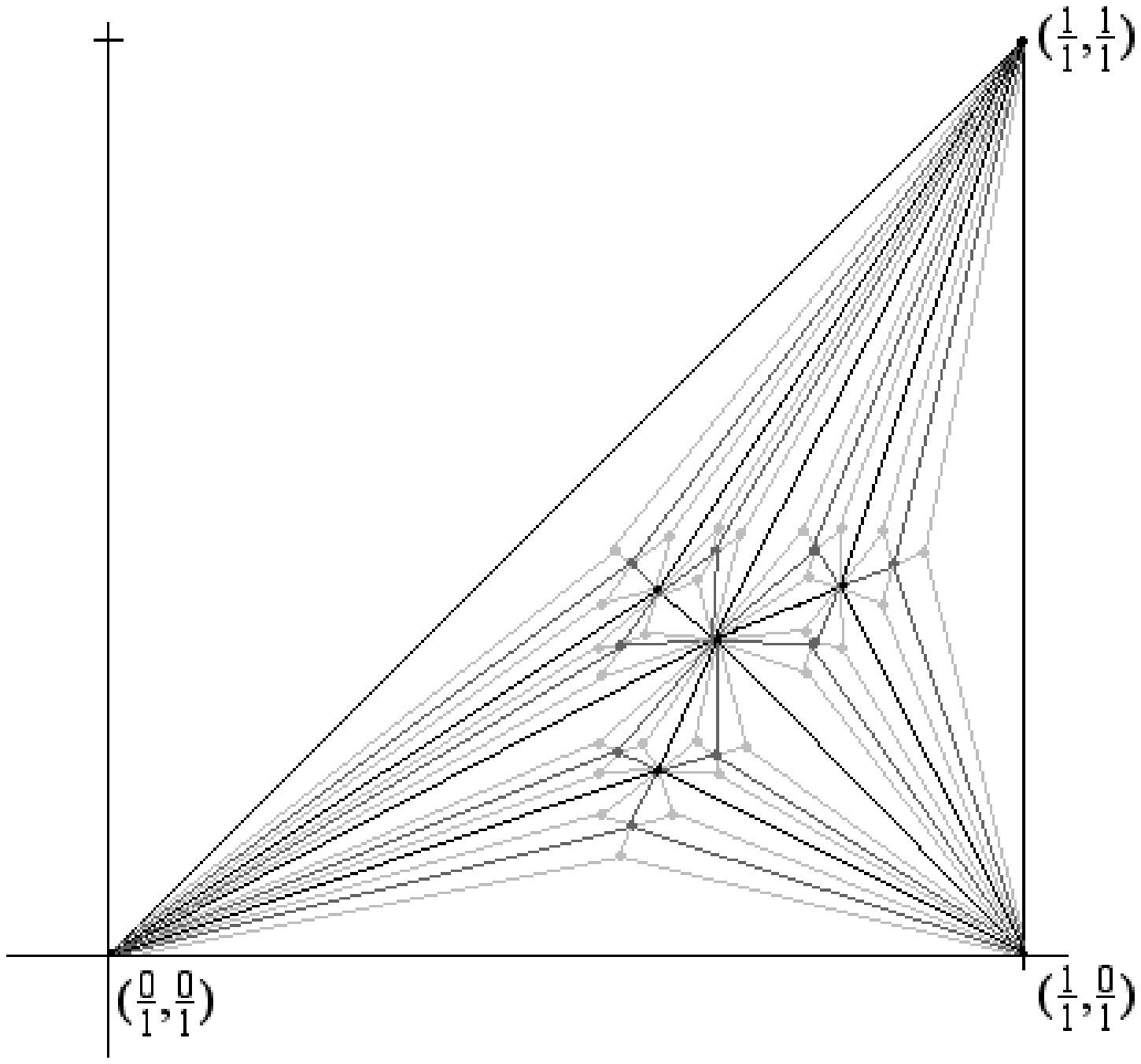}}
\bigskip
\end{tabbing}

\pagebreak

\noindent and for the Bary partition:

\begin{tabbing}
\bigskip
\hspace*{0.375in}\=
\scalebox{0.5}
{\includegraphics*[bb=0in 3.25in 8.5in 9.75in]{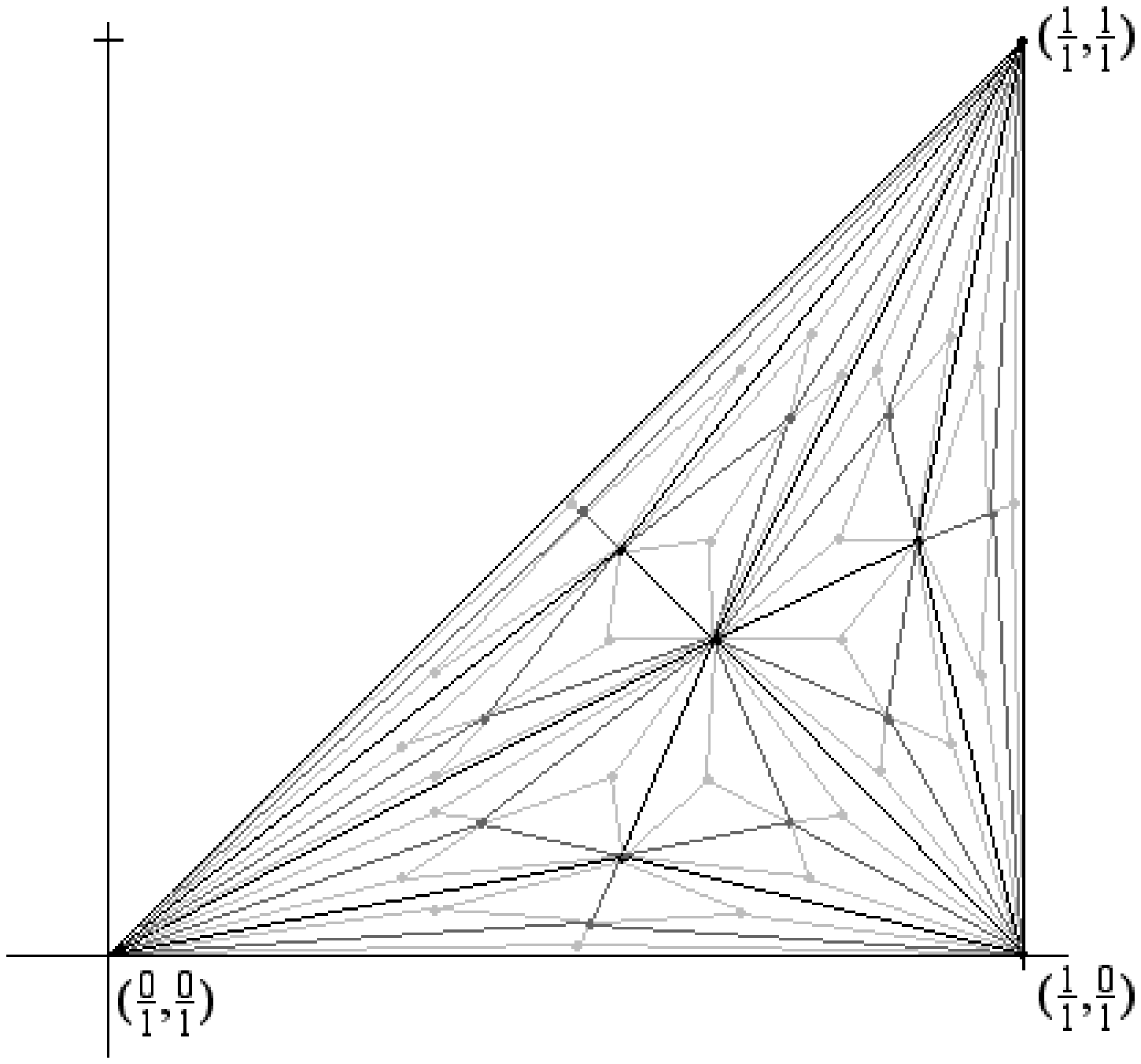}}
\bigskip
\end{tabbing}

\noindent (We find it interesting that the diagram for the Farey 
partition is much more aesthetically pleasing than the one for the 
barycentric partition.)

\noindent By definition, $\delta_{3}(v_{i}(2))=\delta_{2}(v_{i}(2)$ for any 
vertex in a subtriangle $\langle v_{1}(2), v_{2}(2), 
v_{3}(2)\rangle,$ of $\bigtriangleup_{2,F}.$ Thus to describe $\delta_{3},$ we need only specify 
what happens on the new vertices obtained in 
$\bigtriangleup_{3,F}\mbox{ and }\bigtriangleup_{3,B}.$

This new correspondence becomes:
\begin{eqnarray*}
\delta_{3} \pmatrix{11/27\cr 4/27\cr } &=& \pmatrix{5/9\cr 2/9\cr }\\
\delta_{3} \pmatrix{11/27\cr 7/27\cr } &=& \pmatrix{5/9\cr 3/9\cr }\\
\delta_{3} \pmatrix{14/27\cr 1/27\cr } &=& \pmatrix{4/9\cr 1/9\cr }\\
\delta_{3} \pmatrix{14/27\cr 13/27\cr } &=& \pmatrix{4/9\cr 3/9\cr }\\
\delta_{3} \pmatrix{20/27\cr 4/27\cr } &=& \pmatrix{6/9\cr 2/9\cr }\\
\delta_{3} \pmatrix{20/27\cr 16/27\cr } &=& \pmatrix{6/9\cr 4/9\cr }\\
\delta_{3} \pmatrix{23/27\cr 7/27\cr } &=& \pmatrix{7/9\cr 3/9\cr }\\
\delta_{3} \pmatrix{23/27\cr 16/27\cr } &=& \pmatrix{7/9\cr 4/9\cr }\\
\delta_{3} \pmatrix{26/27\cr 13/27\cr } &=& \pmatrix{6/9\cr 3/9\cr }.
\end{eqnarray*}  

\begin{theorem}The sequence of functions ${\{\delta_{n}\}}$ is uniformly 
convergent.
\end{theorem}
\noindent {\bf Proof:}
For any point $v\in \bigtriangleup_{F}$, its image $\delta_{n}(v)$, 
for any $n$, must land in one of the $3^{n}$ subtriangles in the partition $\tilde{\mathcal{P}}_{n},$ the
$n^{th}$ partition of $\bigtriangleup_{B}.$   Label this triangle by $\langle{\tilde{v}_{1}}^{s}(n), {\tilde{v}_{2}}^{s}(n),
{\tilde{v}_{3}}^{s}(n)\rangle$, for $s=1,\ldots, 3^{n}$.  By  the definition, we can see that for any 
$m>n$, the image $\delta_{m}(v)$, while rarely equal to 
$\delta_{n}(v)$, remains in the triangle $\langle{\tilde{v}_{1}}^{s}(n), {\tilde{v}_{2}}^{s}(n),
{\tilde{v}_{3}}^{s}(n)\rangle$

Each subtriangle $\langle{\tilde{v}_{1}}^{s}(n), {\tilde{v}_{2}}^{s}(n),
{\tilde{v}_{3}}^{s}(n)\rangle$, for $s=1,\ldots, 3^{n}$ of the partition 
$\tilde{\mathcal{P}}_{n},$ will  gain a new barycenter in the next step. 
We will denote the barycenter of each such  partitioned triangle
by ${\tilde{v}_{0}}^{s}(n+1),$ where, of course,
$${\tilde{v}_{0}}^{s}(n+1)={\tilde{v}_{1}}^{s}(n)\hat{+}{\tilde{v}_{2}}^{s}(n)
\hat{+}{\tilde{v}_{3}}^{s}(n).$$

Let $\epsilon > 0$ be given.  Clearly, there exist an $N$ such that for the
$N^{th}$ partition of $\bigtriangleup_{B},$ the maximum
distance between any vertex of the $s^{th}$ subtriangle and its new barycenter
is given by
$$\mbox{max}{\{d({\tilde{{v}}_{i}}^{s}(N),{\tilde{{v}}_{0}}^{s}(N+1))}\}\leq\frac{\epsilon}{4},$$
for each $s\in{\{1,2,\dots,3^{N}\}}.$  Also, given any points $u$ and 
$v$ in the subtriangle $\langle
{\tilde{v}_{1}}^{s}(N), {\tilde{v}_{2}}^{s}(N),
{\tilde{v}_{3}}^{s}(N)\rangle,$ we have
\begin{eqnarray*} d(u,w) &\leq&
d(u,{\tilde{v}_{0}}^{s}(N+1))+d({\tilde{v}_{0}}^{s}(N+1),w)\\
&\leq & \mbox{2 max 
d}({\tilde{v}_{i}}^{s}(N),{\tilde{v}_{0}}^{s}(N+1)\\
&\leq &\frac{\epsilon}{2}.
\end{eqnarray*}

Now, let $v\in\bigtriangleup_{F}.$  Then
$\delta_{N}(v)\in\langle{\tilde{v}_{1}}^{s}(N), {\tilde{v}_{2}}^{s}(N),
{\tilde{v}_{3}}^{s}(N)\rangle$ for some $s =1,2,\dots,3^{N}.$  In
particular, for all m,n
$\geq N$, we have $\delta_{m}(v)$ and $\delta_{n}(v)$ also in the 
subtriangle
$\langle{\tilde{v}_{1}}^{s}(N),
{\tilde{v}_{2}}^{s}(N),
{\tilde{v}_{3}}^{s}(N)\rangle.$  But this implies that
$$d(\delta_{m}(v),\delta_{n}(v))\leq\frac{\epsilon}{2}.$$
  Thus
${\{\delta_{n}\}}$ is uniformly Cauchy and the result follows.  $\Box$

\begin{definition} Define the Farey-Bary map to be 
$\delta:\bigtriangleup_{F}\rightarrow\bigtriangleup_{B}$ where 
$\delta$ is the limit of the sequence ${\{\delta_{n}\}}.$
\end{definition}
\begin{theorem}
The Farey-Bary map is continuous.
\end{theorem}
\noindent{\bf Proof:}  Clearly, since each $\delta_{n}$ is a linear map, the 
sequence ${\{\delta_{n}\}}$ consists of continuous functions.  The 
result follows.
$\Box$


\section{Farey Iteration in the Domain as Multi-dimensional Continued
Fraction}
Minkowski's $?(x)$ provides a link  between algebraic 
properties of numbers and the failure of differentiabilty, almost 
everywhere, for $?(x)$.  Our goal is to find analogous links for the 
Farey-Bary
map $\delta(x,y)$.  The key algebraic property of the 
Minkowski question mark function is that $?(x)$ maps
 quadratic irrationals to rational numbers.  The goal of this section 
 is to show that $\delta(x,y)$ maps a class of pairs of cubic 
 irrationals to pairs of rationals.  Unfortunately, we cannot make 
 the claim that $\delta$ maps all pairs of cubics (even in the same 
 number field) to pairs of rationals.

\subsection{Preliminary Notation}
Let $(\alpha, \beta) \in \bigtriangleup_{F}$.  The Farey partitions
of $\bigtriangleup_{F}$ yield a sequence of triangles converging to the
point $(\alpha, \beta)$.  Suppose that at the $n^{th}$ stage of the Farey 
partitioning,
the triangle that contains $(\alpha, \beta)$ is $\langle v_{1}(n), v_{2}(n), 
v_{3}(n)\rangle.$  We will maintain the notation $v_{i}(n)$ to mean either 
the cartesian version of the vertex or 
the vector in space that corresponds 
to the vertex.  That is, $v_{i}(n)$ will refer to
$\pmatrix{p_{i}(n)/r_{i}(n)\cr q_{i}(n)/r_{i}(n)},$ as well as to $\pmatrix{p_{i}(n)\cr q_{i}(n)\cr r_{i}(n)\cr}.$
Furthermore, we will  order the vertices so that for all $n$, 
$$r_{1}(n)\leq r_{2}(n)\leq
r_{3}(n).$$

We want to relate the vertices of the $(n-1)^{st}$ subtriangle that 
contains$(\alpha, \beta)$ with the vertices of the subtriangle at the next
iteration.  For that, suppose that $(\alpha, \beta)\in\langle v_{1}(n-1), v_{2}(n-1), 
v_{3}(n-1)\rangle\subseteq\bigtriangleup_{n-1,F}.$
Applying the next partition, $\mathcal{P}_{n},$ to 
$\bigtriangleup_{n-1,F}$,
we decompose $\langle 
v_{1}(n-1), v_{2}(n-1), v_{3}(n-1)\rangle$ into three new subtriangles.  
If we let $\langle v_{1}(n), v_{2}(n), v_{3}(n)\rangle$ denote the 
subtriangle into which $(\alpha, \beta)$ falls, we see that there are 
three possibilities for the vertices of $\langle v_{1}(n), v_{2}(n), v_{3}(n)\rangle.$  


In case I, the vertices of the newly partitioned triangle will be:
\begin{eqnarray*}
v_{1}(n)& = &v_{1}(n-1) \\
v_{2}(n) &=& v_{2}(n-1) \\
v_{3}(n) &=& v_{1}(n-1)\hat{+}v_{2}(n-1)\hat{+}v_{3}(n-1)
\end{eqnarray*}

\begin{tabbing}
\bigskip
\hspace*{0.3750in}\=
\scalebox{0.5}
{\includegraphics*[bb=0in 4.25in 8.5in 9.75in]{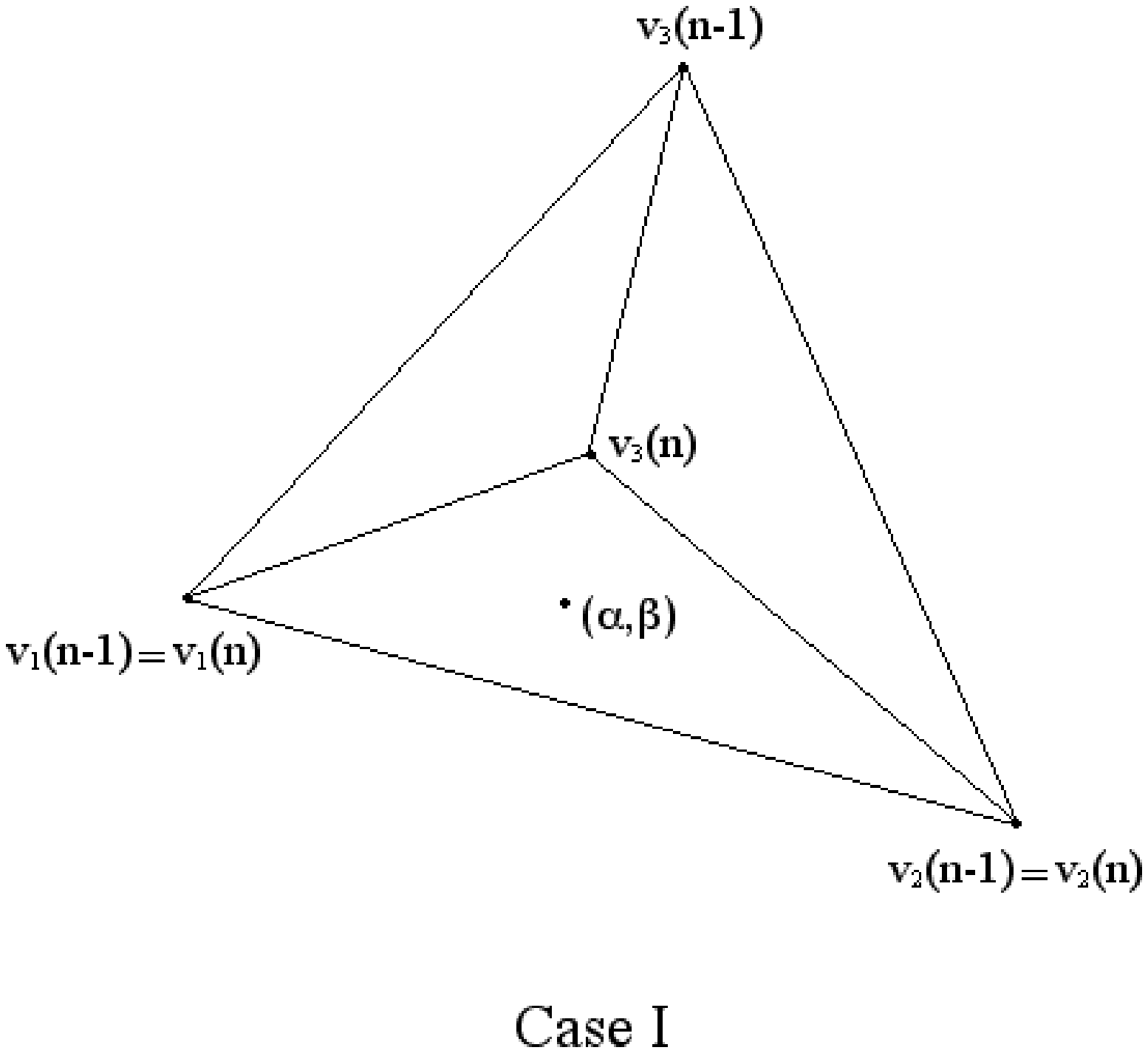}}
\bigskip
\end{tabbing}

Similarly, the vertices in case II will be:
\begin{eqnarray*}
v_{1}(n)& = &v_{2}(n-1) \\
v_{2}(n) &=& v_{3}(n-1) \\
v_{3}(n) &=& v_{1}(n-1)\hat{+}v_{2}(n-1)\hat{+}v_{3}(n-1).
\end{eqnarray*}

\begin{tabbing}
\bigskip
\hspace*{0.3750in}\=
\scalebox{0.5}
{\includegraphics*[bb=0in 4.25in 8.5in 9.75in]{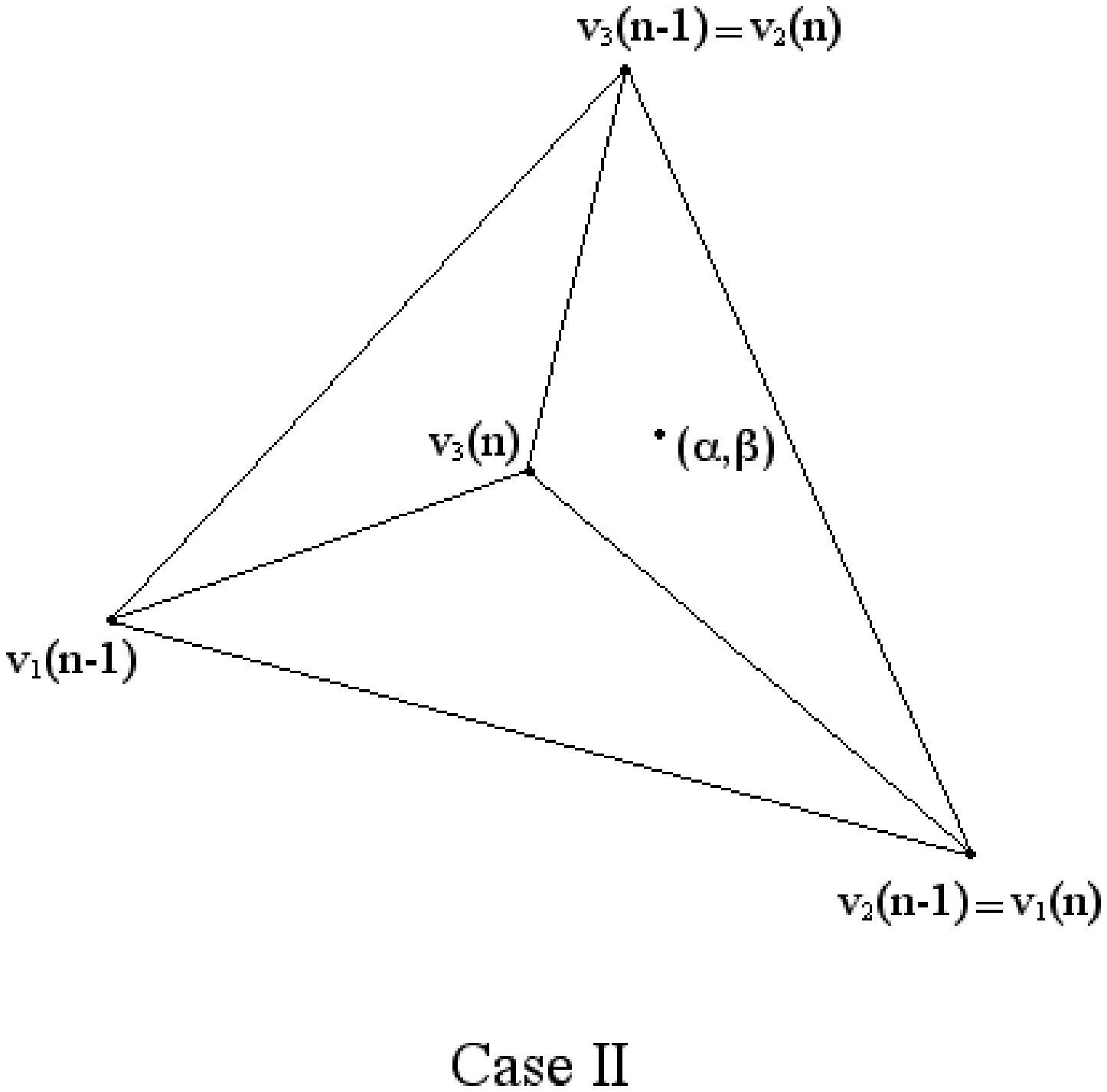}}
\bigskip
\end{tabbing}

For case III we have:
\begin{eqnarray*}
v_{1}(n)& = &v_{1}(n-1) \\
v_{2}(n) &=& v_{3}(n-1) \\
v_{3}(n) &=& v_{1}(n-1)\hat{+}v_{2}(n-1)\hat{+}v_{3}(n-1).
\end{eqnarray*}

\begin{tabbing}
\bigskip
\hspace*{0.3750in}\=
\scalebox{0.5}
{\includegraphics*[bb=0in 4.25in 8.5in 9.75in]{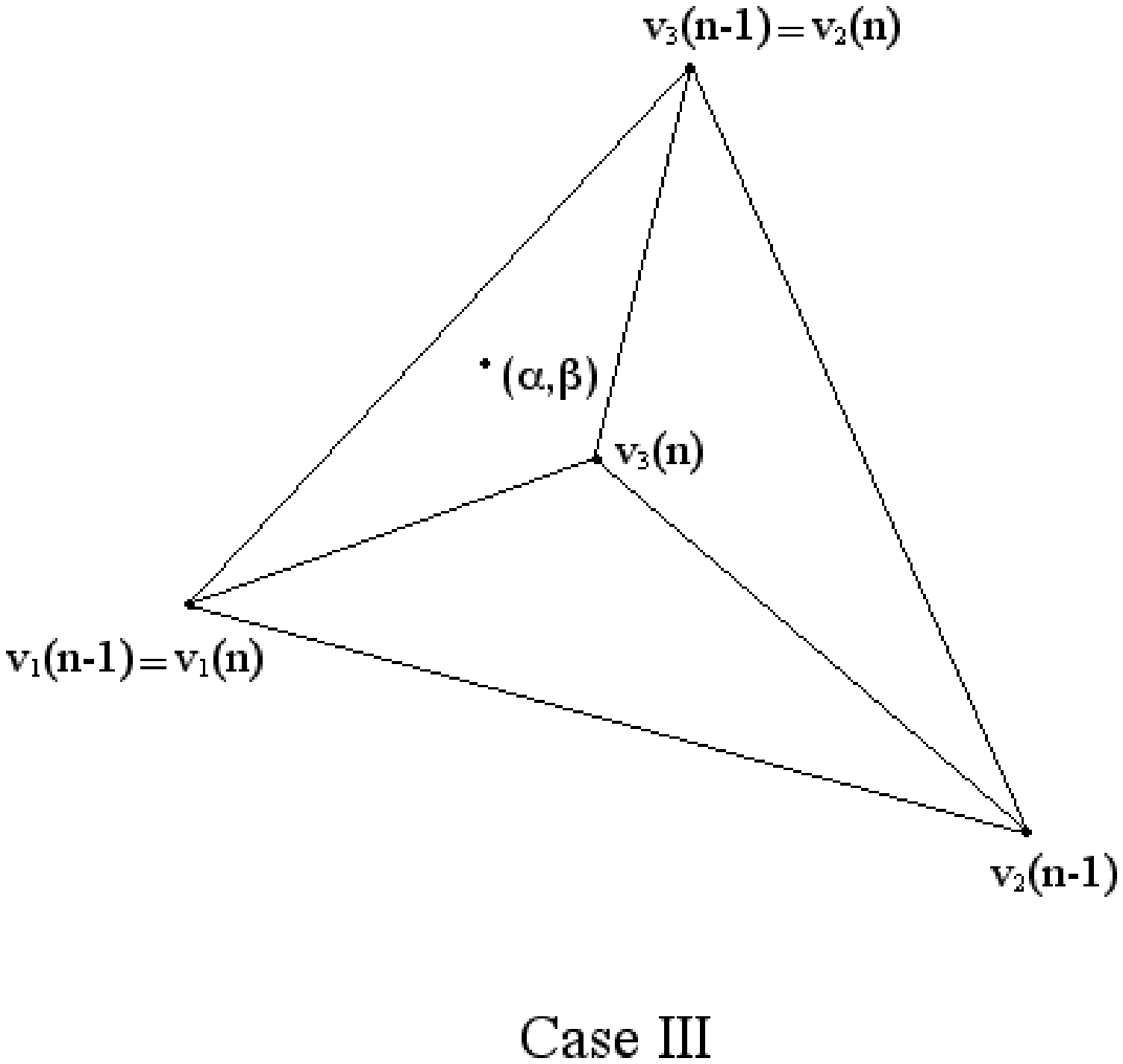}}
\bigskip
\end{tabbing}


In the next section, we will streamline this notation.


\subsection{Fixing Notation}
For each $(\alpha, \beta) \mbox{ in } \bigtriangleup_{F}$ we now 
associate a sequence of positive integers that will uniquely determine the precise 
convergence of the Farey subtriangles to $(\alpha, \beta).$ 

To motivate the eventual notation, consider the following three 
possibilities.  Start with a triangle, with vertices 
$v_{1}$,$v_{2}$, and $v_{3}$, still keeping the convention that 
$r_{1}\leq r_{2} \leq r_{3}$.  Suppose we perform $k$ type I 
operations in a row.  The resulting new triangle will have vertices 
in the following form:
$$v_{1},v_{2}, kv_{1}\hat{+}kv_{2}\hat{+}v_{3}.$$
If we perform a type II operation on the triangle, and then $k-1$ 
type I operations, the new triangle will have vertices:
$$v_{2},v_{3}, v_{1}\hat{+}kv_{2}\hat{+}kv_{3}.$$
If we perform a type III operation on the triangle, and then $k-1$ 
type I operations, the new triangle will have vertices:
$$v_{1},v_{3}, kv_{1}\hat{+}v_{2}\hat{+}kv_{3}.$$
This suggests the following notation.

Define a sequence
$\{a_{1}(i_{1}), a_{2}(i_{2}), \ldots\}$ to be such that each $a_{k}(i_{k})$ 
is a positive integer and each $i_{k}$ represents either case I, II or 
III.  The value of $a_{k}(i_{k})$ denotes the operation of first 
applying a type $i_{k}$ operation and then $a_{k}(i_{k})- 1$ type I 
operations.  We use the further convention that for $k\geq 2$, 
$i_{k}$ can only be of type II or III. 

Note that, in the notation of the previous section, by the time we 
are at step
$a_{k}(i_{k})$, we have performed $n = a_{1}(i_{1})+ a_{2}(i_{2}) 
+\ldots + a_{k}(i_{k})$ Farey partitions of $\bigtriangleup_{F}.$  We associate to each
$(\alpha, \beta) \in \bigtriangleup_{F}$ the sequence that
yields the corresponding Farey partitions that converge to 
$(\alpha,\beta)$.  This sequence will be unique.  We will also use $\langle v_{1}(k), v_{2}(k), 
v_{3}(k)\rangle$ to denote the subtriangle of  $\bigtriangleup_{n,F}$ containing $(\alpha, 
\beta)$ after $n = a_{1}(i_{1})+ a_{2}(i_{2}) 
+\ldots + a_{k}(i_{k})$ steps partitioning $\bigtriangleup_{F}.$  
Finally, if we know what case we are in, that is, if we know $i_{k},$ 
we will simply write $a_{k}$ instead of $a(i_{k}).$

\begin{example}
\end{example}
\noindent The shaded region below corresponds to all points 
$(2(III),1(II), 1(I))$.  Note that in the notation of last section $n=4$ 
but that in the notation of this section and for the rest of the paper $k=3.$


\begin{tabbing}
\bigskip
\hspace*{0.50in}\=
\scalebox{0.5}
{\includegraphics*[bb=0in 3.25in 8.5in 9.75in]{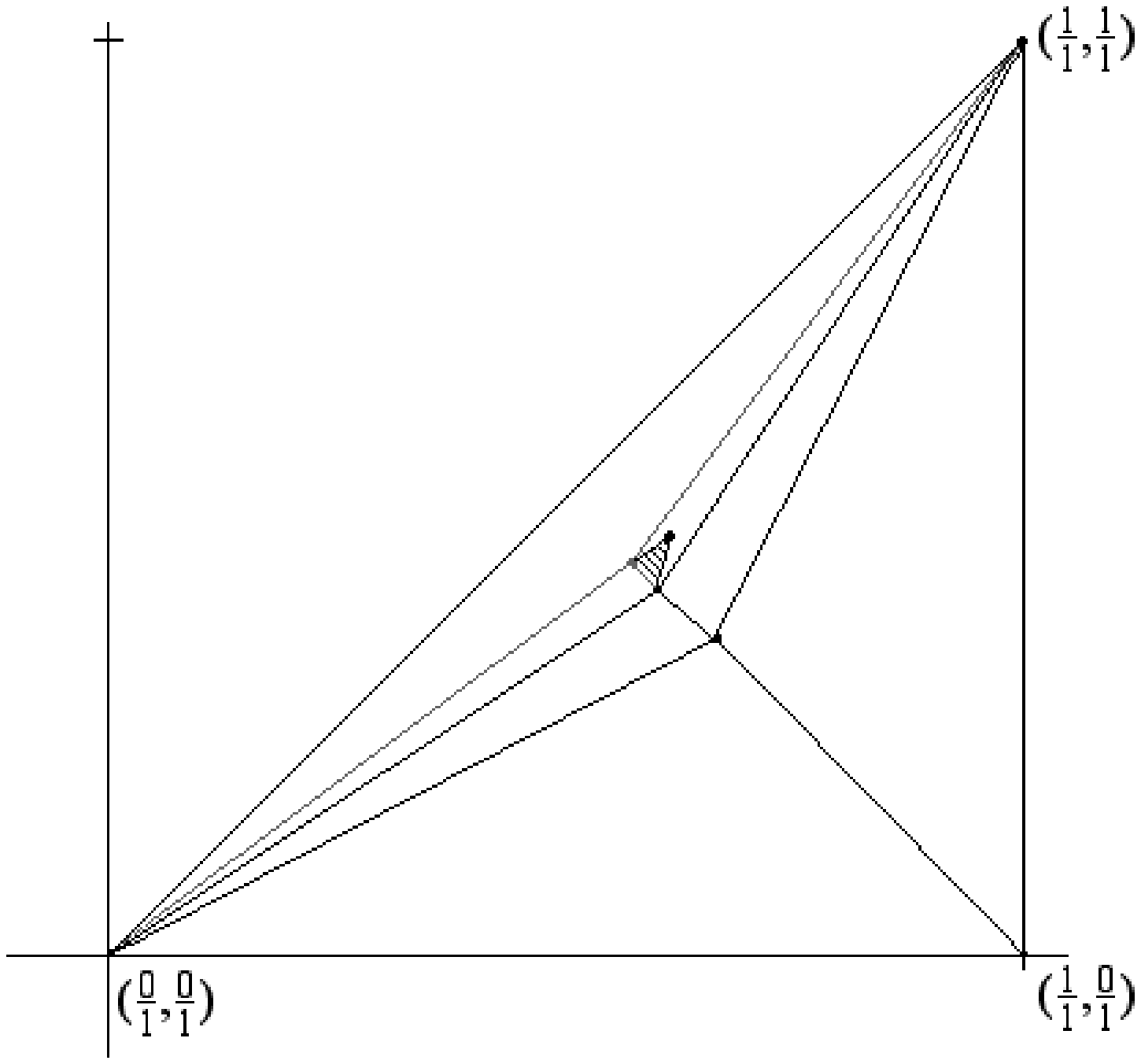}}
\bigskip
\end{tabbing}

We now have the following recursion formulas for the vertices.  For case
I at the k$^{th}$ step we get:
\begin{eqnarray*}
v_{1}(k)& = &v_{1}(k-1) \\
v_{2}(k) &=& v_{2}(k-1) \\
v_{3}(k) &=& a_{k}v_{1}(k-1)\hat{+}a_{k}v_{2}(k-1)\hat{+}v_{3}(k-1),
\end{eqnarray*}
\noindent For case II we have:
\begin{eqnarray*}
v_{1}(k)& = &v_{2}(k-1) \\
v_{2}(k) &=& v_{3}(k-1) \\
v_{3}(k) &=& v_{1}(k-1)\hat{+}a_{k}v_{2}(k-1)\hat{+}a_{k}v_{3}(k-1),
\end{eqnarray*}
Finally, for case III we get:
\begin{eqnarray*}
v_{1}(k)& = &v_{1}(k-1) \\
v_{2}(k) &=& v_{3}(k-1) \\
v_{3}(k) &=& a_{k}v_{1}(k-1)\hat{+}v_{2}(k-1)\hat{+}a_{k}v_{3}(k-1).
\end{eqnarray*}

We can  put these recursion relations naturally into a matrix
language.  At each step $k$, define $M_{k}$ to be the three-by-three
matrix
$$M_{k}=(v_{1}(k)\; v_{2}(k)\;v_{3}(k)) =
\pmatrix{p_{1}(k)&p_{2}(k)&p_{3}(k)\cr
q_{1}(k)&q_{2}(k)&q_{3}(k)\cr
r_{1}(k)&r_{2}(k)&r_{3}(k)}$$
If, from the (k-1)$^{st}$ step to the k$^{th}$ step, we are in case I, then
$$M_{k}=M_{k-1}\pmatrix{1&0&a_{k}\cr
0&1&a_{k}\cr
0&0&1},$$
for case II we have:
$$M_{k}=M_{k-1}\pmatrix{0&0&1\cr
1&0&a_{k}\cr
0&1&a_{k}},$$
and for case III,
$$M_{k}=M_{k-1}\pmatrix{1&0&a_{k}\cr
0&0&1\cr
0&1&a_{k}}.$$
Denote in each of these cases the matrix on the right by
$A_{k}(I),A_{k}(II)$ and $A_{k}(III)$, respectively.
Then we have that each $M_{k}$ is the product of $M_{0}$ with a sequence
of various $A_{m}$.

\begin{theorem}
Each $M_{k}$ is in the special linear group {\bf SL(3,Z)}.
\end{theorem}

\noindent {\bf Proof:} All we need to show is that for all $k$,
$$\det(M_{k})= \pm 1.$$
This follows immediately from observing that $\det(M_{0})= 1$ and
that the determinants of each of the various $A_{k}(I), A_{k}(II)$
and $A_{k}(III)$ are also plus or minus one.  $\Box$


\subsection{Areas of Farey Subtriangles}
Given a finite sequence $\{a_{1}(i_{1}), a_{2}(i_{2}), \ldots,
a_{k}(i_{k})\}$ of positive integers, we
define
$$\bigtriangleup_{k}=\{(x,y):\{a_{1}(i_{1}), \ldots, a_{k}(i_{k})\} 
\;\mbox{are the}\; 1^{\mbox{st}}\;k\;\mbox{terms in Farey sequence}\}.$$

A major goal of this paper is showing that the
areas of these triangles $\bigtriangleup_{k }$ cannot go to zero too 
quickly.  For these
calculations, we will need an easy formula for
the areas of the $\bigtriangleup_{k}$.

\begin{theorem}
The area of a triangle with vertices $(p_{1}/r_{1}, q_{1}/r_{1})$,
$(p_{2}/r_{2}, q_{2}/r_{2})$, and  $(p_{3}/r_{3}, q_{3}/r_{3})$ is
$$\mbox{Area of triangle}= \frac{1}{2}
\frac{|\det \pmatrix{p_{1}&p_{2}&p_{3}\cr
q_{1}&q_{2}&q_{3}\cr
r_{1}&r_{2}&r_{3}}|}{r_{1}r_{2}r_{3}}.$$
\end{theorem}
This is just a calculation involving cross products.
\begin{cor}
Given any finite sequence $\{a_{1}(i_{1}), a_{2}(i_{2}), \ldots,
a_{k}(i_{k})\}$ of positive integers,
$$2\mbox{Area}\;(\bigtriangleup_{k})=
\frac{1}{r_{1}(k)r_{2}(k)r_{3}(k)}.$$
\end{cor}
This follows since $\det(M_{k})= \pm 1.$


\subsection{Farey Periodicity Implies Cubic Irrationals}

As we iterate our procedure, the vertices of our
triangles converge to a single vector.
We want to show:

\begin{theorem}
Suppose that $(\alpha, \beta)\in \bigtriangleup_{F}$ has an eventually
periodic Farey sequence.  Then both $\alpha$ and $\beta$ are
algebraic numbers with $deg(\alpha)\leq3,$ $deg(\beta)\leq3$ and
$$\dim_{{\bf Q}}{\bf Q}[\alpha, \beta]\leq 3.$$
\end{theorem}
This is why the Farey partitioning can be viewed as a 
multi-dimensional continued fraction algorithm.

\noindent {\bf Proof:}
We will be  heavily  using two facts.  First,  an eigenvector $(1,a,b)$ of a $3\times 3$
matrix with rational coefficients has the property that
$$\dim_{{\bf Q}}{\bf Q}[a,b]\leq 3,$$ as seen in a similar argument 
in \cite{Garrity1} in section eight.
Second,  if we multiply a matrix, which has a  largest real
eigenvalue, repeatedly by itself, in the limit the columns of
the matrix converge to the eigenvector corresponding to the largest
eigenvalue.

Suppose that $(\alpha, \beta)\in \bigtriangleup$ has an eventually
periodic Farey sequence.  Even if it is not periodic, the vertices of
the corresponding Farey partition triangles converge to the point $(\alpha,
\beta)$.
We have seen above that the vertices of the partition triangles
correspond to the columns of  matrices that are the products of
various $A_{k}(I)$, $A_{k}(II)$ and $A_{k}(III)$.  With the assumption
of periodicity, denote the product of the initial non-periodic
matrices be $B$ and the product of the periodic part be $A$.  Then
some of the
Farey partition triangles about the point $(\alpha, \beta)$ are given
by
$$B,BA,BA^{2},BA^{3},\ldots.$$
The columns of the matrices $A,A^{2},A^{3},\ldots$
 must converge to a multiple of $B^{-1}(1,\alpha,
\beta)^{T}$. But the columns of the $A^{k}$ must also converge to an 
eigenvector and  hence   $B^{-1}(1,\alpha,
\beta)^{T}$ is an eigenvector of the matrix $A$.  This will give 
us that
 $\alpha$ and $ \beta$ must have the desired properties.   $\Box$


\section{Iteration in the Barycentric Range}

We have defined Farey partitions, $\mathcal{P}_{n},$ in 
$\bigtriangleup_{F}$ and barycentric partitions, $\tilde{\mathcal{P}}_{n}$ 
in $\bigtriangleup_{B}.$  In $\bigtriangleup_{F},$ the Farey partitions 
yielded an association between each $(\alpha,\beta)$ and a sequence 
obtained from the convergence of the subtriangles resulting from 
successive applications of the partitions, $\mathcal{P}_{n}.$  This 
association depended only on the successive partitioning of each subtriangle 
into three more subtriangles and not on the relative positioning of 
the new subtriangles.  We can follow the same procedure in 
$\bigtriangleup_{B}.$  That is, if we let $(a,b)\in \bigtriangleup_{B},$ we can 
again associate with $(a,b)$ a sequence of positive integers
$\{\tilde{a}_{1}(i_{1}), \tilde{a}_{2}(i_{2}), \ldots\}$ which 
come from a sequence of barycentric triangles converging to the point $(a,b).$

Label the triangle corresponding to
$\{\tilde{a}_{1}(i_{1}), \tilde{a}_{2}(i_{2}), \ldots, 
\tilde{a}_{k}(i_{k}\}$ by
$$\tilde{\bigtriangleup}\{\tilde{a}_{1}(i_{1}), \tilde{a}_{2}(i_{2}), \ldots,
\tilde{a}_{k}(i_{k}\}.$$  Recall that $$\tilde{\bigtriangleup}_{B}=\langle\tilde{v}_{1}(0), \tilde{v}_{2}(0), 
\tilde{v}_{3}(0)\rangle,$$ where $$\tilde{v}_{1}(0)= \pmatrix{0\cr 0\cr 
1\cr}, \tilde{v}_{2}(0)=\pmatrix{1\cr
0\cr 1\cr}\mbox{ and }\tilde{v}_{3}(0)=\pmatrix{1\cr 1\cr 1\cr}.$$

Associated with the sequence $\{\tilde{a}_{1}(i_{1}), \tilde{a}_{2}(i_{2}), 
\ldots, \tilde{a}_{k}(i_{n}\}$
will be vertices ${v_{1}}(k),
{v_{2}}(k)$ and ${v_{3}}(k)$ and corresponding vectors
$$v_{1}(k)= \pmatrix{*\cr *\cr 3^{a_{1}+\dots +a_{k}}\cr},
v_{2}(k)=\pmatrix{*\cr
*\cr 3^{a_{1}+\dots +a_{k}}\cr},
 v_{3}(k)=\pmatrix{*\cr *\cr 3^{a_{1}+\dots +a_{k}}\cr},$$
 where the other entries for the vectors are nonnegative integers.
 There are, of course, matrices $\tilde{M_{n}}$ that map the vertices from a
given
 level to the vertices of the next level, in analogue to the matrices
$M_{n}$.  The $\tilde{M_{n}}$ are products of matrices of the form:
$$\pmatrix{3&0&1\cr
0&3&1\cr
0&0&1}, \pmatrix{0&0&1\cr
3&1&1\cr
0&3&1}, \pmatrix{3&0&1\cr
0&0&1\cr
0&3&1}$$
depending if we are in case I, II or III, respectively.

Note that at each individual step of the barycentric partitioning, we are
cutting the area down by a factor of a $1/3$.  This leads to the 
following theorem.


\begin{theorem}
Twice the area of $\tilde{\bigtriangleup}\{\tilde{a}_{1}(i_{1}), \tilde{a}_{2}(i_{2}), \ldots,
 \tilde{a}_{k}(i_{k})\}$ is $\frac{1}{3^{\tilde{a}_{1} + \tilde{a}_{2} + \dots +
 \tilde{a}_{k}(i_{k})}}$.
 \end{theorem}

 
\subsection{Ternary Periodicity implies rationality}

Suppose that we have a point $(a,b) \in \bigtriangleup_{B}$ for which the
barycentric partitioning is eventually periodic.  We want to show
that both $a$ and $b$ are rational numbers. That is, we want the following 
theorem.

\begin{theorem}
If $(a,b)\in \bigtriangleup_{B}$ has an eventually
perdiodic Barycentric sequence, then both $a$ and $b$ are
rational.
\end{theorem}

\noindent {\bf Proof:}
This proof is almost exactly the same as the
corresponding proof for the Farey case, whose notation we adopt.  There is one significant
difference, namely that the matrices  whose columns yield the
vertices of the barycentric partitioning are all multiples of a
stochastic matrix.  This means that each matrix is a multiple of a
matrix whose columns add to one.  If the columns add to one, then it
can easily be shown that the limit of the products of such a matrix
converges to a matrix whose rows are multiples of $(1,1,1)$ (see 
chapter six in \cite{Minc1}). Thus the matrices $A, A^{2}, 
A^{3},\ldots$ converge to a matrix whose rows are multiples of 
$(1,1,1)$.  Since everything in sight is rational, we can show that 
$B^{-1}(1,\alpha,\beta)^{T}$ will converge to a triple of rational 
numbers.   Since the entries 
of $B$ are integers, this yields that $\alpha$ and $\beta$ are 
rational numbers.  $\Box$


\section{The Farey-Bary Analog of Singularness}


The original Minkowksi ?(x) function is singular, meaning that even 
though it is increasing and continuous, it has derivative zero almost 
everywhere.  The key to the proof lies in showing that at almost all 
points
$$\liminf_{k\rightarrow \infty}\frac{\mbox{length of interval in 
range}}{\mbox{length of interval in domain}}=0,$$
for appropriately defined intervals.  We will show a direct analog of 
this, where the lengths of intervals are replaced by areas of 
triangles.  Thus we will show 
$$\liminf_{k\rightarrow \infty}\frac{\mbox{area of subtriangle in 
range}}{\mbox{area of subtriangle  in domain}}=0,$$
again for appropriately defined subtriangles.


\subsection{Almost everywhere $\limsup \frac{a_{1}+\ldots a_{n}}{n}=\infty$}
This is the most technically difficult section of the paper.  The goal
is to show the following theorem, which will be critical in the next 
section.  Recall that given any point $(\alpha, \beta)\in 
\bigtriangleup$, we have associated a sequence 
$\{a_{1},a_{2},\ldots\}$ of positive integers.  We want to show that this sequence must 
increase to infinity, in some sense, almost everywhere.  The precise 
statement is:
\begin{theorem}
The set of $(\alpha, \beta)\in \bigtriangleup$ for which
$$\limsup_{n\rightarrow \infty} \frac{a_{1}+\dots +a_{n}}{n}<
\infty$$
has measure zero.
\end{theorem}
As it will only be apparent in the next section  why to we need this theorem, we 
recommend on the first reading of this paper to go 
to the next section first.  

Before proving the theorem, we need a preliminary lemma.
First, let
$$v_{1}= \pmatrix{x_{1}\cr
y_{1}\cr
z_{1}},v_{2}=\pmatrix{x_{2}\cr
y_{2}\cr
z_{2}}, v_{3}=\pmatrix{x_{3}\cr
y_{3}\cr
z_{3}}$$
and let 
$$T=\langle v_{1}, v_{2}, v_{3}\rangle$$
 be the corresponding triangle in the plane, where 
the $v_{i}$ are now viewed as points in the plane.

Suppose that $\det{v_{1} v_{2} v_{3}}=1$.  Then we know that
$$2\cdot \mbox{area of} \langle v_{1}, v_{2}, 
v_{3}\rangle =\frac{1}{z_{1}z_{2}z_{3}}.$$
Given a positive number $L>1$, define
$T_{L}(1)$ to be the triangle with vertices $\hat{v}_{1}, \hat{v}_{2}$
and $Lv_{1}\hat{+}Lv_{2}\hat{+}v_{3}$, $T_{L}(2)$
the triangle with vertices $\hat{v}_{2},
v_{1}\hat{+}Lv_{2}\hat{+}Lv_{3}$
and $\hat{v}_{3}$  and $T_{L}(3)$
the triangle with vertices
$\hat{v}_{1}, Lv_{1}\hat{+}v_{2}\hat{+}Lv_{3}$
and $\hat{v}_{3}$.

\begin{tabbing}
\bigskip
\hspace*{0.675in}\=
\scalebox{0.5}
{\includegraphics*[bb=0in 5.in 8.5in 9.75in]{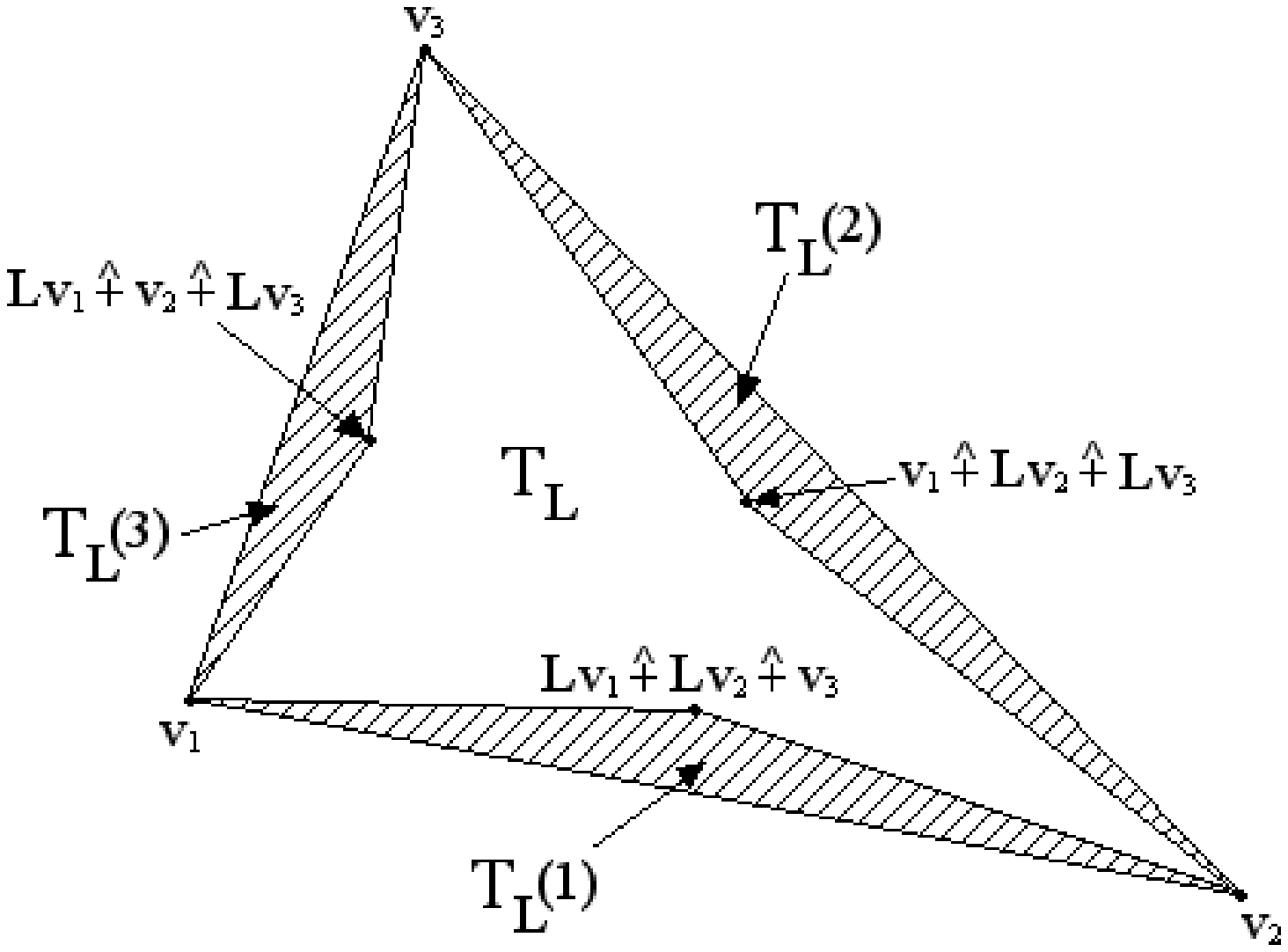}}
\bigskip
\end{tabbing}

\noindent Define
$$T_{L}=T-T_{L}(1)-T_{L}(2)-T_{L}(3).$$
We now state and then prove a lemma that is the technical heart of
the proof of the
theorem:
\begin{lemma}For all $L\geq 1$,
$$\mbox{area}(T_{L})\leq \frac{L-1}{L}\mbox{area}(T).$$
\end{lemma}

\noindent {\bf Proof of Lemma:}
We know that
$$2\cdot \mbox{area}(T) = \frac{1}{z_{1}z_{2}z_{3}}.$$
For ease of notation, we set $z_{1}=x,z_{2}=y, z_{3}=z$.
Then
$$2\cdot \mbox{area}(T_{L}) =
\frac{1}{xyz}-\frac{1}{x(Lx+y+Lz)z}-\frac{1}{xy(Lx+Ly+z)}-\frac{1}{(x+Ly+Lz)yz}$$
$$= \frac{1}{xyz}[1- \frac{y}{(Lx+y+Lz)}-\frac{z}{(Lx+Ly+z)}
-\frac{x}{(x+Ly+Lz)}].$$
Thus we must show that
$$[1- \frac{y}{(Lx+y+Lz)}-\frac{z}{(Lx+Ly+z)}
-\frac{x}{(x+Ly+Lz)}]\leq \frac{L-1}{L}.$$
Setting
$$\alpha = x^{2}y+x^{2}z+xy^{2}+xz^{2}+y^{2}z+yz^{2},$$
we have that
$$[1- \frac{y}{(Lx+y+Lz)}-\frac{z}{(Lx+Ly+z)}
-\frac{x}{(x+Ly+Lz)}]$$
$$=\frac{\alpha (L^{3}-L) + 2xyz (L^{3}-1)}{(Lx +
Ly +z)(Lx + y + Lz)(x+Ly+Lz)}.$$
After a series of calculations, we get that this is equal to:
$$(L-1)[\frac{L(L+1)\alpha + 2xyz(L^{2}+L+1)}{L(L^{2}+L+1)\alpha + 
L^{2}(x^{3}+y^{3}+z^{3}) +
xyz(3L^{2}+2L^{3}+1)}].$$
Thus we must show that
$$[\frac{L(L+1)\alpha + 2xyz(L^{2}+L+1)}{L(L^{2}+L+1)\alpha +
L^{2}(x^{3}+y^{3}+z^{3})+
xyz(3L^{2}+2L^{3}+1)}]\leq \frac{1}{L},$$
which is equivalent to showing that
$$L^{2}(L+1)\alpha + 2xyzL(L^{2}+L+1) \leq $$
$$L(L^{2}+L+1)\alpha +
L^{2}(x^{3}+y^{3}+z^{3}) + xyz(2L^{3}+3L^{2}+ 1),$$
which in turn, reduces to showing that
$$2Lxyz\leq \alpha + L^{2}(x^{3}+y^{3}+z^{3}) + L^{2}xyz + xyz.$$
This last inequality follows from the fact that $L\geq 1$.  Thus the proof of the
lemma is done.

\noindent {\bf Proof of Theorem:}
For each positive integer $N$, set
$$M_{N}= \{(\alpha,\beta)\in \bigtriangleup:\mbox{for all}\; n\geq 1,
\frac{a_{1}+\dots + a_{n}}{n}\leq N\}.$$
We will show that
$$\mbox{measure}(M_{N})=0.$$
Since the union of all of the $M_{N}$ is the set we want to show has
measure zero, we will be done.

Now,
$\frac{a_{1}+\dots + a_{n}}{n}\leq N$ if and only if
$$a_{1}+\dots+ a_{n}\leq nN.$$
Since each $a_{i}\geq 1$, this last inequality implies
 $$n-1 + a_{n}\leq nN$$
 or
 $$a_{n}\leq n(N-1)+1.$$
 Set
 $$\tilde{M}_{N}=\{(\alpha,\beta)\in \bigtriangleup:\mbox{for all}\; n\geq 1,
 a_{n}\leq n(N-1)+1\}.$$
 Since $M_{N}\subset \tilde{M}_{N}$, if  we can show that
 measure$(\tilde{M})_{N}=0$, we will be done.

 Set
 $$\tilde{M}_{N}(1)=\{(\alpha,\beta)\in \bigtriangleup:
 a_{1}\leq (N-1)+1\}$$
 and in general
 $$\tilde{M}_{N}(k)=\{(\alpha,\beta)\in \tilde{M}_{N}(k-1):
 a_{k}\leq k(N-1)+1\}$$
 Then we have a decreasing nested sequence of sets with
 $$\tilde{M}_{N} = \bigcap_{k=1}^{\infty}\tilde{M}_{N}(k).$$
 But this puts us into the
 language of the above lemma.  Letting $L=k(N-1)+1$,  we can conclude that

 $$\mbox{measure}(\tilde{M}_{N}(k))\leq 
 \frac{k(N-1)}{k(N-1)+1}\mbox{measure}(\tilde{M}_{N}(k-1))$$
 and hence
 
 $$\mbox{measure}(\tilde{M}_{N})\leq \prod_{k=2}^{\infty}\frac{k(N-1)}{k(N-1)+1}.$$
 We must show
 this infinite product is zero, which is equivalent to showing that 
 its reciprical
 $$\prod_{k=2}^{\infty}\frac{k(N-1)+1}{k(N-1)}=
 \prod_{k=2}^{\infty}(1+\frac{1}{k(N-1)})=\infty.$$
 Taking logarithms, this is the same as showing that the series
 $$\sum_{k=2}^{\infty}\log(1+\frac{1}{k(N-1)})=\infty.$$
 This in turn follows since, for large enough $k$, we have
 $$\log(1+\frac{1}{k(N-1)})\geq \frac{1}{2k(N-1)}.$$
 We are done.


\subsection{Almost everywhere $$\liminf (\mbox{area}(\tilde{\Delta}_{n})/ 
\mbox{area}(\Delta_{n})) = 0$$}
The goal of this section, and for the entire paper, is:
\begin{theorem}
For  any point $(\alpha, \beta)\in \bigtriangleup$, off of a set of
measure zero,
$$\liminf_{n\rightarrow \infty} 
\frac{\mbox{area}(\tilde{\bigtriangleup}\{a_{1}(i_{1}),
a_{2}(i_{2}), \ldots,
 a_{n}(i_{n})\})}{\mbox{area}(\bigtriangleup\{a_{1}(i_{1}), a_{2}(i_{2}), \ldots,
 a_{n}(i_{n})\})} = 0.$$
 \end{theorem}
 This is capturing the intuition that the determinant of the Jacobian
 of the map $\delta:\bigtriangleup_{F} \rightarrow \bigtriangleup_{B}$ is
 zero almost everywhere, which in turn is a direct generalization that
 the Minkowski question-mark function is singular.  In fact,
  our proof is in spirit a generalization of Viader,
Paradis and Bibiloni's  work in \cite{Viader-Paradis-Bibiloni1}.

 \noindent {\bf Proof:} We know that, letting $s_{n}=a_{1}+\dots
 a_{n},$
 $$\mbox{area}(\tilde{\bigtriangleup}\{a_{1}(i_{1}), a_{2}(i_{2}), \ldots,
 a_{n}(i_{n})\}) = \frac{1}{2 \cdot 3^{s_{n}}}$$
 and that
 $$\mbox{area}(\bigtriangleup\{a_{1}(i_{1}), a_{2}(i_{2}), \ldots,
 a_{n}(i_{n})\}) = \frac{1}{2\cdot r_{1}(n) r_{2}(n)r_{3}(n)}.$$
 Thus we want to show that, almost everywhere,
 $$\liminf_{n\rightarrow \infty}
 \frac{r_{1}(n)r_{2}(n)r_{3}(n)}{3^{s_{n}}} = 0.$$

  We know  that $r_{3}(n)$ is
  $$a_{n}r_{1}(n-1) + a_{n}r_{2}(n-1) + r_{3}(n-1),$$
  $$a_{n}r_{1}(n-1) + r_{2}(n-1) + a_{n}r_{3}(n-1),$$
  or
  $$r_{1}(n-1) + a_{n}r_{2}(n-1) + a_{n}r_{3}(n-1).$$
  Thus we have, 
  by the convention of our notation, $$r_{1}(n)\leq r_{2}(n) \leq 
 r_{3}(n)\leq (2a_{n}
 +1)r_{3}(n-1).$$
 By iterating this inequality, we have
 $$r_{1}(n)\leq r_{2}(n) \leq r_{3}(n)\leq \prod_{i=1}^{n}(2a_{j}+1).$$
 Thus
$$\frac{\mbox{area}(\tilde{\bigtriangleup}\{a_{1}(i_{1}),
a_{2}(i_{2}), \ldots,
 a_{n}(i_{n})\})}{\mbox{area}(\bigtriangleup\{a_{1}(i_{1}), a_{2}(i_{2}), \ldots,
 a_{n}(i_{n})\})}\leq \frac{\prod_{i=1}^{n}(2a_{j}+1)^{3}}{3^{s_{n}}}.$$

By the arithmetic-geometric mean, 
 $$\prod_{i=1}^{n}(1+b_{i})\leq (1 + \frac{b_{1}+\dots
 b_{n}}{n})^{n}.$$
 Setting $b_{j}=2a_{j}$, we get
 \begin{eqnarray*}
 \frac{\mbox{area}(\tilde{\bigtriangleup}\{a_{1}(i_{1}),
a_{2}(i_{2}), \ldots,
 a_{n}(i_{n})\})}{\mbox{area}(\bigtriangleup\{a_{1}(i_{1}), a_{2}(i_{2}), \ldots,
 a_{n}(i_{n})\})}&\leq& \frac{(1+\frac{2s_{n}}{n})^{3n}}{3^{s_{n}}} \\
 &\leq& \frac{(\frac{3s_{n}}{n})^{3n}}{3^{s_{n}}} \\
 &\leq& (\frac{27\cdot (\frac{s_{n}}{n})^{3}}{3^{s_{n}/n}})^{n}.
 \end{eqnarray*}

From the previous section, we know that $s_{n}/n \rightarrow \infty$,
almost everywhere.  Since the above denominator has a $3^{s_{n}/n}$ 
term while the numerator only has a $(s_{n}/n)^{3}$ term, the entire 
ratio must approach zero, giving us our result.

\section{Questions}

There are a number of natural questions.  First, all of this can 
almost certainly be generalized to higher dimensions.  

More importantly, how much does the function theory  of
$\delta$ influence the diophantine properties of points 
in $\bigtriangleup$?   

There are many multi-dimensional continued fraction algorithms.  For 
any of these that involve partitioning a given triangle into three new 
subtriangles, a map analogous to our $\delta$ can of course be 
defined.  What are the properties of these new maps?

Underlying most work on multidimensional continued fractions, though 
frequently hidden behind view, are Lie theoretic properties of the 
special linear group.      Can this be made more explicit?

Finally, the initial Hermite problem remains open.

\end{document}